\documentclass{amsart}
\usepackage{amssymb,latexsym,enumitem,mathrsfs}
\usepackage{hyperref}
\usepackage{color}
\usepackage{stmaryrd}
\usepackage{lineno}
\usepackage{bbm}   
\usepackage{scalerel,stackengine}
\usepackage{tikz}

\def\lesim{\lesssim}
\def\beq{\begin{equation}}
\def\eeq{\end{equation}}

\newcommand\avsuminner[2]{%
  {\sbox0{$\m@th#1\sum$}%
   \vphantom{\usebox0}%
   \ooalign{%
     \hidewidth
     \smash{\vrule height\dimexpr\ht0+1pt\relax depth\dimexpr\dp0+1pt\relax}%
     \hidewidth\cr
     $\m@th#1\sum$\cr
   }%
  }%
}

\makeatletter

\newcommand*\barredprod{%
  \DOTSB\mathop{%
      \@rodriguez@mathpalette \@rodriguez@overprint@bar \prod
    }\slimits@
}
\newcommand*\@rodriguez@mathpalette[2]{%
  \mathchoice
    {#1\displaystyle      \textfont         {#2}}%
    {#1\textstyle         \textfont         {#2}}%
    {#1\scriptstyle       \scriptfont       {#2}}%
    {#1\scriptscriptstyle \scriptscriptfont {#2}}%
}
\newcommand*\@rodriguez@overprint@bar[3]{%
  \sbox\z@{$#1#3$}%
  \dimen@   = \ht\z@   \advance \dimen@   \p@
  \dimen@ii = \dp\z@   \advance \dimen@ii \p@
  \dimen4 = 1.25\fontdimen 8 #2\thr@@ \relax
  \ooalign{% the resulting box has the same...
    \@rodriguez@bar \dimen@ \z@ \cr   % ... height as the first row
    $\m@th #1#3$\cr
    \@rodriguez@bar \z@ \dimen@ii \cr % ... depth as the last row
  }%
}

\newcommand*\@rodriguez@bar[2]{%
  \hidewidth \vrule \@width \dimen4 \@height #1\@depth #2\hidewidth
}

\makeatother
\def\pnorm#1{ \Big(  #1 \Big) }

\def\anorm#1{ \Big|  #1 \Big| }
\def\Norm#1{ \Big\|  #1 \Big\| }
\def\norm#1{\big\|  #1 \big\| }
\def\inn#1#2{\langle#1,#2\rangle}
\def\set#1{ \Big\{ #1 \Big\} }
\def\dim{\textup{dim}}

\def\sgn{\mathrm{sgn}}

\newcommand{\pisharp}{\pi_{\theta\#}}
\newcommand{\aaa}{a}
\newcommand{\heavy}{\mathrm{heavy}}
\newcommand{\light}{\mathrm{light}}
\usepackage{esint}% See geometry.pdf to learn the layout options. There are lots.
\usepackage{graphicx}
\newtheorem{theorem}{Theorem}
\newtheorem{lemma}{Lemma}
\newtheorem{remark}{Remark}
\newtheorem{claim}{Claim}

\newcommand{\cC}{\mathcal{C}}
\newcommand{\cD}{\mathcal{D}}

\newcommand{\cH}{\mathcal{H}}

\newcommand{\R}{\mathbb R}

\newcommand{\C}{\mathbb C}
\newcommand{\N}{\mathbb N}
\newcommand{\T}{\mathbb T}

\renewcommand{\S}{\mathbb S}

\newcommand{\W}{\mathbb{W}}
\newcommand{\e}{\varepsilon}

\newcommand{\s}{\psi}

\renewcommand{\d}{\delta}

\newcommand{\mc}{\mathcal}

\def\set4{\mathcal I}
\def\tup14{(1,2,3,4)}

\def\be{{\mathbf e}}

\newcommand\vwidehat[1]{\arraycolsep=0pt\relax%
\begin{array}{c}
\stretchto{
  \scaleto{
    \scalerel*[\widthof{\ensuremath{#1}}]{\kern-.5pt\bigwedge\kern-.5pt}
    {\rule[-\textheight/2]{1ex}{\textheight}} %WIDTH-LIMITED BIG WEDGE
  }{\textheight} % 
}{0.5ex}\\           % THIS SQUEEZES THE WEDGE TO 0.5ex HEIGHT
#1\\                 % THIS STACKS THE WEDGE ATOP THE ARGUMENT
\rule{-1ex}{0ex}
\end{array}
}
\newtheorem*{comm*}{Comment}
\newtheorem{definition}{Definition}
\newtheorem*{lemma*}{Lemma}
\newtheorem{corollary}{Corollary}

\newtheorem{proposition}{Proposition}

\usepackage{bbm}   
\usepackage{scalerel,stackengine}
\stackMath
\newcommand\widecheck[1]{%
\savestack{\tmpbox}{\stretchto{%
  \scaleto{%
    \scalerel*[\widthof{\ensuremath{#1}}]{\kern-.6pt\bigwedge\kern-.6pt}%
    {\rule[-\textheight/2]{1ex}{\textheight}}%WIDTH-LIMITED BIG WEDGE
  }{\textheight}% 
}{0.5ex}}%
\stackon[1pt]{#1}{\scalebox{-1}{\tmpbox}}%
}
\usepackage{ esint }

\newcommand{\supp}{\mathrm{supp}}

\newcommand{\ZZ}{\mathbb{Z}}

\newcommand{\D}{\mathbb D}
\newcommand{\de}{\delta}
\newcommand{\ga}{\gamma}
\newcommand{\Ga}{\Gamma}
\newcommand{\ZS}{\mathbb S}

\newcommand{\wh}{\widehat}
\newcommand{\wt}{\widetilde}
\newcommand{\si}{\sigma}
\newcommand{\Si}{\Sigma}
\newcommand{\Tau}{\mathcal T}

\newcommand{\en}{\epsilon_{\circ}}

\newcommand{\dist}{\textup{dist}}

\newcommand{\Id}{\boldsymbol 1}

\newcommand{\U}{\mathbb {U}}

\begin{document}

\author[Gan]{Shengwen Gan}
\address{Department of Mathematics\\
Massachusetts Institute of Technology\\
Cambridge, MA 02142-4307, USA}
\email{shengwen@mit.edu}

\author[Guo]{Shaoming Guo}
\address{Department of Mathematics\\
University of Wisconsin-Madison\\
Madison, WI-53706, USA}
\email{shaomingguo@math.wisc.edu}

\author[Guth]{Larry Guth}
\address{Department of Mathematics\\
Massachusetts Institute of Technology\\
Cambridge, MA 02142-4307, USA}
\email{lguth@math.mit.edu}

\author[Harris]{Terence L.~J.~Harris}
\address{Department of Mathematics\\ Cornell University\\ Ithaca\\ NY 14853-4201\\ USA}
\curraddr{Department of Mathematics\\
University of Wisconsin-Madison\\
Madison, WI-53706, USA}
\email{terry.harris@wisc.edu}

\author[Maldague]{Dominique Maldague}
\address{Department of Mathematics\\
Massachusetts Institute of Technology\\
Cambridge, MA 02142-4307, USA}
\email{dmal@mit.edu}

\author[Wang]{Hong Wang}
\address{Department of Mathematics\\
UCLA\\
Los Angeles, CA 90095, USA}
\email{hongwang@math.ucla.edu}

\keywords{decoupling inequalities, superlevel set}
\subjclass[2020]{42B15, 42B20}

\date{}

\title{On restricted projections to planes in $\R^3$ }
\maketitle

\begin{abstract}
    Let $\gamma:[0,1]\rightarrow \mathbb{S}^{2}$ be a non-degenerate curve in $\mathbb{R}^3$, that is to say, $\det\big(\gamma(\theta),\gamma'(\theta),\gamma''(\theta)\big)\neq 0$. For each $\theta\in[0,1]$, let
    $V_\theta=\gamma(\theta)^\perp$ and let
    $\pi_\theta:\mathbb{R}^3\rightarrow V_\theta$ be the orthogonal projections. We prove that if $A\subset \mathbb{R}^3$ is a Borel set, then for a.e. $\theta\in [0,1]$ we have $\textup{dim}(\pi_\theta(A))=\min\{2,\textup{dim} A\}$.
    
    More generally, we prove an exceptional set estimate. For $A\subset\R^3$ and $0\le s\le 2$, define $E_s(A):=\{\theta\in[0,1]: \textup{dim}(\pi_\theta(A))<s\}$. We have $\textup{dim}(E_s(A))\le \max\{1+s-\textup{dim}(A),0\}$.
    We also prove that if $\textup{dim}(A)>2$, then for a.e. $\theta\in[0,1]$ we have $\mathcal{H}^2(\pi_\theta (A))>0$.
\end{abstract}

\tableofcontents

\section{Introduction}

Let $\ZS^2\subset \R^3$ denote the unit sphere. Let $\gamma: [0, 1]\to\ZS^2$ be a $C^2$ curve. We say that $\gamma$ is non-degenerate if 
\begin{equation}
    \det(\gamma(\theta), \gamma'(\theta), \gamma''(\theta))\neq 0,
\end{equation}
for every $\theta\in [0, 1]$. A model example for the non-degenerate curve is $\ga_\circ:\theta\mapsto (\frac{\cos\theta}{\sqrt{2}},\frac{\sin\theta}{\sqrt{2}},\frac{1}{\sqrt{2}})$ $(\theta\in[0,1])$.

For a given $\theta\in [0, 1]$, let $V_{\theta}=\gamma(\theta)^{\perp}$ denote the orthogonal complement of $\gamma(\theta)$ in $\R^3$, and let $\pi_{\theta}: \R^3\to \gamma(\theta)^{\perp}$ denote the orthogonal projection onto $\gamma(\theta)^{\perp}$. For $\alpha>0$ and a Borel set $E\subset \R^3$, we will use $\mathcal{H}^{\alpha}(E)$ to denote the $\alpha$-dimensional Hausdorff measure of $E$. Moreover, we use $\dim X$ to denote the Hausdorff dimension of a set $X$.

\begin{theorem}\label{mainplane}
Suppose $A\subset \R^3$ is a Borel set of Hausdorff dimension $\alpha$. For $0\le s< 2$, define the exceptional set 
\[ E_s=\{\theta\in[0,1]: \dim(\pi_\theta(A))<s \}. \]
Then we have
$$ \dim (E_s)\le \max\{1+s-\alpha,0\}.  $$
\end{theorem}

As an immediate corollary, we have:

\begin{corollary}
Suppose $A\subset \R^3$ is a Borel set of Hausdorff dimension $\alpha$. Then we have
$$ \dim(\pi_\theta(A))=\min\{2,\alpha\},\textup{~for~a.e.~}\theta\in[0,1]. $$
\end{corollary}

\begin{theorem}\label{thmpositive}
Suppose $A\subset \R^3$ is a Borel set of Hausdorff dimension greater than $2$. Then  
\begin{equation}
    \mathcal{H}^2(\pi_{\theta}(A))>0, 
\end{equation}
for almost every $\theta\in [0, 1]$. 
\end{theorem}

% If $\ga:[0,1]\rightarrow\ZS^2$ is a curve that satisfies the non-degenerate condition $$\det\big(\ga(\theta),\ga'(\theta),\ga''(\theta)\big)\neq 0,$$ 
% then 
% we call $\ga$ a \textit{non-degenerate curve}.
% A model example for the non-degenerate curve is
% $\ga_\circ:\theta\mapsto (\frac{\cos\theta}{\sqrt{2}},\frac{\sin\theta}{\sqrt{2}},\frac{1}{\sqrt{2}})$ $(\theta\in[0,1])$. 

% In this paper, we study the the projections in $\R^3$ whose directions are determined by $\ga$.
% For each $\theta\in[0,1]$, let $V_\theta\subset \R^3$ be the $2$-dimensional subspace that is orthogonal to $\ga(\theta)$.
% We also define $\pi_\theta:\R^3\rightarrow V_\theta$ to be the orthogonal projection onto $V_\theta$. 
% We use $\dim X$ to denote the Hausdorff dimension of set $X$. Let us state our main results.

\subsection{Background of the problems}\hfill\\
The projection theory dates back to Marstrand \cite{marstrand1954some}, who showed that if $A$ is a Borel set in $\R^2$, then the projection of $A$ onto almost every line through the origin has Hausdorff dimension $\min\{1,\dim A\}$. This was generalized to higher dimensions by Mattila \cite{mattila1975hausdorff}, who showed that if $A$ is a Borel set in $\R^n$, then the projection of $A$ onto almost every $k$-plane through the origin has Hausdorff dimension $\min\{k,\dim A\}$.
It is more general to consider the projection problem when the direction set is restricted to some submanifold of the Grassmannian. To have a better understanding of this restricted projection problem, the first step is to study the problem in $\R^3$. 
F\"assler and Orponen made a conjecture about restricted projections to lines and planes (see Conjecture 1.6 in \cite{fassler2014restricted}), and there has been much related research (see for example \cite{fassler2014restricted}, \cite{chen2018restricted}, \cite{jarvenpaa2008one}, \cite{jarvenpaa2014hausdorff}, \cite{kaenmaki2017marstrand}, \cite{oberlin2015application}, \cite{orponen2015hausdorff}, \cite{orponen2020improved}, \cite{harris2019improved}, \cite{harris2021restricted}). For more of an introduction to this problem, we refer to \cite{harris2019improved}. Recently, Käenmäki-Orponen-Venieri \cite{kaenmaki2017marstrand} and Pramanik-Yang-Zahl \cite{pramanik2022furstenberg} proved one half of the conjecture: restricted projections to lines. In this paper, we resolve another half of the conjecture: restricted projections to planes.

\subsection{An overview of the high-low method}
The high-low method is a powerful tool developed recently in Fourier analysis. There are many applications of the high-low method, see for example \cite{guth2019incidence}, \cite{demeter2020small}, \cite{guth2020sharp}, \cite{guth2020improved}, \cite{fu2021sharp}, \cite{gan2022square}.

In this subsection, we briefly discuss how the high-low method can be used to study projection theory. As a warm-up, we study Marstrand's projection theorem from another point of view, using the high-low method.

\begin{theorem}[Marstrand's projection theorem]
For each $\theta\in[0,\pi]$, define $L_\theta:=\{x\in\R^2: \arg(x)=\theta\}$ and let $p_\theta:\R^2\rightarrow L_\theta$ be the projection.
Suppose $A\subset \R^2$ is a Borel set, then we have $\dim (p_\theta(A))=\min\{1,\dim A\}$ for a.e. $\theta\in[0,\pi]$.
\end{theorem}

We will frequently use the following definition.
\begin{definition}
For a number $\de>0$ and any set $X$, we use $|X|_\de$ to denote the maximal number of $\de$-separated points in $X$.
\end{definition}

Marstrand's projection theorem can be reduced to the following discretized version. We do not show how reduction works in this section, but we will give the full details in the later section when we prove our main theorems.
\begin{proposition}\label{prop}
Fix $0<s<1$.
Fix a small scale $\de>0$ and set $\Theta=\de\ZZ\cap[0,\pi]$. For each $\theta\in\Theta$, define $L_\theta:=\{x\in\R^2: \arg(x)=\theta\}$ and let $p_\theta:\R^2\rightarrow L_\theta$ be the projection.
Suppose $A\subset B^2(0,1)$ is a union of disjoint  $\de$-balls with measure $|A|=\de^{2-a}$, or equivalently $|A|_\de\sim \de^{-a}$. We also assume there is a subset $\Theta'\subset \Theta$ with $\#\Theta'\gtrsim \de^{-1}$, such that for any $\theta\in\Theta'$, $p_\theta(A)$ (which is a union of line segments of length $\de$ in $L_\theta$) satisfies the $s$-dimensional condition: 
     For each $r\ge \de$ and line segment $I_r\subset L_\theta$ of length $r$, we have 
     \begin{equation}\label{condi2}
         |p_\theta(A)\cap I_r|_\de\lesssim   (r/\de)^s.
     \end{equation}
Then,
\begin{equation}
    \de^{-a}\lesssim_{s} \de^{-s}.
\end{equation}
\end{proposition}

% \begin{remark}
% {\rm
% We explain how it is related to Marstrand's projection theorem. By discretization, the set $A$ is like an $a$-dimensional set at scale $\de$. We only need to consider $a\le 1$. Our goal is to show for most of $\theta\in\Theta$, $p_\theta(A)$ is at least an ``$a$-dimensional" set. Suppose by contradiction, there is an $s<a$ such that $p_\theta(A)$ satisfies the s-dimensional condition as in Proposition \ref{prop} for $\theta$ in a large subset $\Theta'\subset \Theta$. By Proposition \ref{prop}, we have $\de^{-a}\lesssim \de^{-s}$, which contradicts $s<a$ as $\de\rightarrow 0$. The  Marstrand's projection theorem in $\R^2$ can be reduced to Proposition \ref{prop} by a discretization argument, which we omit here.
% }
% \end{remark}

\begin{proof}
For each $\theta\in \Theta'$, let $\T_\theta$ be a set of $\de\times 1$ tubes that cover $p_{\theta}^{-1}(p_\theta(A))\cap B^2(0,1)$ and hence cover $A$. We can also assume that $\T_\theta$ satisfies a similar $s$-dimensional condition that is inherited from $p_\theta(A)$.
For each $\theta$, let $S_\theta$ be a $\de^{-1}\times 1$-tube centered at the origin whose longest direction forms an angle $\theta$ with $x$-axis. We see that $S_\theta$ is dual to the tubes in $\T_\theta$.
Now, for each $T_\theta\in\T_\theta,$ we choose a bump function $\psi_\theta$ satisfying the following properties: $\psi_{T_\theta}\ge 1$ on $T_\theta$, $\psi_{T_\theta}$ decays rapidly outside $T_\theta$, and $\supp\wh\psi_{T_\theta}\subset S_\theta$.

Define functions
$$ f_\theta:=\sum_{T_\theta\in\T_\theta}\psi_{T_\theta}\ \ \ \ \textup{and}\ \ \ \ f=\sum_{\theta\in\Theta'} f_\theta. $$
By definition, $f(x)\gtrsim \#\Theta'$ for any $x\in A$, so we simply have
$$ |A|(\#\Theta')^2\lesssim \int_A |f|^2. $$
We can do better by performing a high-low decomposition for $f$.

Let $K$ be a large number that will be determined later. Let $\eta_{low}(\xi)$ be a smooth bump function adapted to $B^2(0,K^{-1}\de^{-1})$ and let $\eta_{high}=1-\eta_{low}$. We have the following high-low decomposition for $f$:
$$ f=f_{low}+f_{high}, $$
where $\wh f_{low}=\eta_{low}\wh f $ and $\wh f_{high}=\eta_{high}\wh f$.
We will show that the high part dominates on $A$. Actually, for $x\in A$, we have
\begin{equation}\label{easyhighlow}
    \de^{-1}\lesssim f(x)\le |f_{high}(x)|+|f_{low}(x)|.
\end{equation}
By definition, $f_{low}(x)=\eta_{low}^\vee * f(x)=\eta_{low}^\vee*\big(\sum_{\theta\in\Theta'}\sum_{T_\theta\in\T_\theta}\psi_{T_\theta}\big)(x)$.
Note that $|\eta_{low}^\vee|$ is morally an $L^1$-normalized bump function at $B^2(0,K\de)$. By the $s$-dimensional condition of $\T_\theta$, we have $$\eta_{low}^\vee*\big(\sum_{T_\theta\in\T_\theta}\psi_{T_\theta}\big)(x)\lesssim K^{-1}\#\{T_\theta\in\T_\theta: T_\theta\cap B^2(x,CK\de)\neq \emptyset\}\lesssim K^{s-1}.$$
In the last step, we use the condition \eqref{condi2} in Proposition \ref{prop}. 
As a result, we have
$$ |f_{low}(x)|\lesssim \#\Theta'K^{s-1}\lesssim \de^{-1}K^{s-1}.$$
Since $s<1$, by choosing $K$ large enough (depending on $s$) and plugging into \eqref{easyhighlow}, we see that
$$ |f_{high}(x)|\ge C^{-1}\de^{-1}-|f_{low}(x)|\gtrsim \de^{-1}. $$
We obtain that
\begin{equation}
    |A|\de^{-2}\lesssim \int_A |f_{high}|^2.
\end{equation}
Next, we will use a strong separation property for the high part. Note that $\wh f_{high}=\sum_{\theta\in\Theta'} \eta_{high} \wh f_\theta$, and  $\{\supp(\eta_{high}\wh f_\theta)\}_\theta$ is at most $O(K)$-overlapped. We have
$$ \int |f_{high}|^2=\int |\wh f_{high}|^2\lesssim K^2 \sum_{\theta\in\Theta'}\int|\eta_{high}\wh f_\theta|^2\lesssim K^2 \sum_{\theta\in\Theta'} \int|f_\theta|^2\lesssim K^2 \de^{-s}.  $$
In the last inequality, we used the $s$-dimensional condition of $p_\theta(A)$.

As a result, we have
$$ |A|\de^{-2}\lesssim_{s} \de^{-s}, $$
which implies $\de^{-a}\lesssim_{s}\de^{-s}$.
\end{proof}

\noindent {\bf Notation.} 
\begin{enumerate}
    \item For two positive real numbers $R_1$ and $R_2$, we say that $R_1\lesssim R_2$ if there exists a large constant $C$, depending on relevant parameters, such that $R_1\le C R_2$; we say that $R_1\ll R_2$ if $R_1\le R_2/C$. 
    \item We use $\dim(E)$ for the Hausdorff dimension of $E$. 
    \item For a given Borel measure $\mu$ supported on $\R^3$ and the projection $\pi_{\theta}$, the pushforward measure $\pisharp\mu$, supported on $\gamma(\theta)^{\perp}$, is defined by 
    \begin{equation*}
        (\pisharp\mu)(E):=\mu((\pi_{\theta})^{-1}(E)),
    \end{equation*}
    for every Borel $E\subset \gamma(\theta)^{\perp}$. 
    \item Let $\mu$ be a compactly supported Borel measure on $\R^3$. Take $\alpha>0$. Define 
\begin{equation*}
    c_{\alpha}(\mu):=\sup_{x\in \R^3, r>0}\frac{\mu(B(x, r))}{r^{\alpha}},
\end{equation*}
where $B(x, r)$ is the ball of radius $r$ centered at $x\in \R^3$. 
\item For $r\ge 1$ and a rectangular box $T\subset \R^3$, we use $rT$ to mean the dilation of $T$ by $r$ with respect to the center of $T$, unless otherwise stated. For $r>0$ and $E\subset \R^3$, we use $r\cdot E$ to mean $\{rx: x\in E\}$. 
% \item The relation among various parameters
% \begin{equation}
%     \delta\ll \delta_0\ll \epsilon
% \end{equation}
\item We often use $m(E)$ the Lebesgue measure of the set $E\subset \R^3$; that is, $m(E)=\mathcal{H}^3(E)$.

\item By a measure we always assume that it is non-negative, unless stated otherwise. 
\end{enumerate}

\bigskip

\begin{sloppypar}
\noindent {\bf Acknowledgement.} S. Guo is partly supported by NSF-1800274 and NSF-2044828. L. Guth is supported by a Simons Investigator Award. H. Wang is supported by NSF Grant DMS-2055544. D. Maldague is supported by the NSF under Award No. 2103249. \end{sloppypar}

\section{Proof of Theorem~\ref{mainplane}}\label{sectionplane}
In this section we prove Theorem \ref{mainplane}.
First, we make some remarks on $\ga$. 
We can cut $\gamma$ into several small pieces and work on each of them. From now on, we assume that $\gamma: [0, \aaa]\to\ZS^2$ is $C^2$ and non-degenerate, and satisfies
\begin{equation}
    \gamma(0)=(0, 0, 1), \ \ \gamma'(0)=(1, 0, 0), \ |\gamma'(\theta)|=1, \forall \theta\in [0, \aaa]. 
\end{equation}
Here $\aaa>0$ is sufficiently small depending on $\gamma$.
Since the parameter $a$ does not play any role here, we may pretend $a=1$. We need some notation.
\begin{definition}[$(\de,s)$-set]\label{deltaset}
Let $P\subset \R^n$ be a bounded set. Let $\de>0$ be a dyadic number, and let $0\le s\le d$. We say that $P$ is a $(\de,s)$-set if
$$ |P\cap B_r|_\de\lesssim  (r/\de)^s ,  $$
for any $B_r$ being a ball of radius $r$ with $\de\le r\le 1$.
\end{definition}

Let $\cH^t_\infty$ denote the $t$-dimensional Hausdorff content which is defined as
$$ \cH^t_\infty(B):=\inf\{ \sum_i r(B_i)^t: B\subset \cup_i B_i \}. $$
Here, each $B_i$ in the covering is a cube and $r(B_i)$ is the length of the cube.
We recall the following  result (see \cite{fassler2014restricted} Lemma 3.13). 
\begin{lemma}\label{lemdelta}
Let $\de,s>0$, and $B\subset \R^n$ with $\cH_\infty^s(B):=\kappa>0$. Then there exists a $\de$-separated $(\de,s)$-set $P\subset B$ with cardinality $\#P\gtrsim \kappa \de^{-s}$.
\end{lemma}

Our main effort will be devoted to the proof of the following theorem.

\begin{theorem}\label{planediscrt} Fix $0<s<2$. For each $\e>0$, there exists $C_{s,\e}$ so that the following holds. Let $\d>0$. Let $H\subset B^3(0,1)$ be a union of  disjoint $\de$-balls and we use $\#H$ to denote the number of $\de$-balls in $H$. Let $\Theta$ be a $\de$-separated subset of $[0,1]$ such that $\Theta$ is a $(\de,t)$-set and $\#\Theta\gtrsim (\log\de^{-1})^{-2}\de^{-t}$ for some $t>0$. Assume for each $\theta\in \Theta$, we have a collection of $\de\times \de\times 1$-tubes $\T_\theta$ pointing in direction $\ga(\theta)$. Each $\T_\theta$ satisfies the $s$-dimensional condition: 
\begin{enumerate}
    \item $\#\T_\theta\lesssim \de^{-s}$, 
    \item $\#\{T\in\T_\theta:T\cap B_r\}\lesssim (\frac{r}{\de})^{s}$, for any $B_r$ being a ball of radius $r$ $(\de\le r\le 1)$.
\end{enumerate}
We also assume that each $\de$-ball contained in $H$ intersects $\gtrsim (\log\de^{-1})^{-2}\#\Theta$ many tubes from $\cup_{\theta\in\Theta} \T_\theta$.
Then 
\[ \#\Theta\#H\le C_{s,\e}\de^{-1-s-\e}.\] 
\end{theorem}

% \begin{theorem}\label{linediscrt} Fix $0<s<1$. For each $\e>0$, there exists $C_{s,\e}$ so that the following holds. Let $\d>0$. Let $A\subset B^3(0,1)$ have size $|A|\ge \d^{-a}\d^3$. Let $\Theta$ be a collection of $\d$-separated $\theta\in[0,1]$ which satisfies: (1) $\#\Theta\gtrsim \de^{-1}$; (2) $|\rho_\theta(A)|_\d\le C_1\d^{-s}$; (3) For each $r\ge \d$ and segment $I_r\subset l_\theta$ of length $r$, we have $|\pi_\theta(A)\cap B_r|_\d\le C_2(r/\d)^{s}$. Then 
% \[ \d^{-a}\le C_{\e,C_1,C_2}\de^{-s-\e}.\] 
% \end{theorem}

\subsection{\texorpdfstring{$\d$}{Lg}-discretization of the projection problem}

In this subsection we show how Theorem \ref{planediscrt} implies Theorem \ref{mainplane}.
Before starting the proof, we state a very useful lemma. We use the following notation. For any $\de=2^{-k}$ ($k\in \mathbb N^+$), let $\cD_\de$ denote the lattice of $\de$-squares in $[0,1]^2$.
For technical reasons, we remove the top edge and the right edge of each $\de$-square so that they are disjoint. 

\begin{lemma}\label{usefullemma}
Suppose $X\subset [0,1]^2$ with $\dim X< s$. Then for any $\e>0$, there exist dyadic squares $\cC_{2^{-k}}\subset \cD_{2^{-k}}$ $(k>0)$ so that 
\begin{enumerate}
    \item $X\subset \bigcup_{k>0} \bigcup_{D\in\cC_{2^{-k}}}D, $
    \item $\sum_{k>0}\sum_{D\in\cC_{2^{-k}}}r(D)^s\le \e$,
    \item $\cC_{2^{-k}}$ satisfies the $s$-dimensional condition: For $l<k$ and any $D\in \cD_{2^{-l}}$, we have $\#\{D'\in\cC_{2^{-k}}: D'\subset D\}\le 2^{(k-l)s}$.
\end{enumerate}
\end{lemma}

\begin{proof}[Proof of the lemma]
Consider all the covering $\cC$ of $X$ by dyadic lattice squares that satisfy condition (1), (2) in Lemma \ref{usefullemma}, i.e., $\cC\subset \bigcup_{k>0}\cD_{2^{-k}}$, $X\subset\bigcup_{D\in\cC}D$ and $\sum_{D\in\cC}r(D)^s\le \e$. We also assume all the dyadic squares in $\cC$ are disjoint. 
We will define an order ``$<$" between any two of such coverings $\cC, \cC'$. First, we define the $k$-th covering number of $\cC$ by
$$ c_k(\cC):=\#(\cC\cap \cD_{2^{-k}}), $$
which is the number of $2^{-k}$-squares in the covering $\cC$.

We say $\cC<\cC'$, if they satisfy: (1) There is a maximal $k_0\ge 0$ such that
$\cC\cap\cD_{2^{-k}}=\cC'\cap\cD_{2^{-k}}$ ($k<k_0$), and $\cC\cap \cD_{2^{-k_0}}\subset \cC'\cap \cD_{2^{-k_0}}$; (2) For any $x\in X$, the square in $\cC'$ that covers $x$ contains the square in $\cC$ that covers $x$. 
It is not hard to check the transitivity: If $\cC<\cC'$ and $\cC'<\cC''$, then $\cC<\cC''$.

Suppose $\cC$ is a covering that is maximal with respect to the order $<$. Then we can show that $\cC$ satisfies condition $(3)$ in Lemma \ref{usefullemma}.
Suppose by contradiction, there exist $l<k$ and $D\in\cD_{2^{-l}}$ so that
\begin{equation}\label{bythis}
    \#\{D'\in \cC\cap\cD_{2^{-k}}:D'\subset D \}>2^{(k-l)s}. 
\end{equation} 
We define another covering $\cC'$ by adding $D$ to $\cC$ and deleting %$\{D'\in\cC\cap\cD_{2^{-k}}:D'\subset D\}$ from $\cC$%.
$\{D'\in\cC \setminus \{D\} :D'\subset D\}$ from $\cC$. It is easy to check that $\cC'$ is still a covering of $X$. By \eqref{bythis}, we can also check $\sum_{D\in\cC'}r(D)^s<\sum_{D\in\cC}r(D)^s\le \e$, so $\cC'$ satisfies $(2)$ in Lemma \ref{usefullemma}. However, $\cC<\cC'$ which contradicts the maximality of $\cC$. 

Now, it suffices to find a maximal element among all the coverings that satisfy condition $(1), (2)$ in Lemma \ref{usefullemma}. First of all, such covering exists by the definition of Hausdorff dimension and $\dim X<s$. By Zorn's lemma, it suffices to find an upper bound for any ascending chain.

Let $\{ \mathcal{C}_j\}_{j \in J}$ be an infinite chain of coverings of $X$. Define 
\[ \mathcal{C} = \bigcap_{j \in J} \bigcup_{\substack{ i \in J \\ \mathcal{C}_i \geq \mathcal{C}_j } } \mathcal{C}_i. \]
We show that $\cC$ is an upper bound of the chain. First, we show that $\cC$ covers $X$.
For $x \in X$ and $j$, let $D^{(j)}_x$ be the largest dyadic square in $\bigcup_{\substack{ i \in J \\ \mathcal{C}_i \geq \mathcal{C}_j } } \mathcal{C}_i$ containing $x$. By the definition of the partial order and the fact that chains are totally ordered, $D^{(j)}_x = D_x$ is independent of $j$, and thus $D_x \in \mathcal{C}$. This shows that $\mathcal{C}$ is a covering of $X$. It also shows that the squares in $\mathcal{C}$ are disjoint. Let $K \in \mathbb{N}$. Choose $j \in J$ such that $\mathcal{C}_i \cap \mathcal{D}_{2^{-k}} = \mathcal{C}_j \cap \mathcal{D}_{2^{-k}}$ for all $0 \leq k \leq K$ and all $\mathcal{C}_i \geq \mathcal{C}_j$. Then 
\[ \sum_{k=0}^K \sum_{D \in \mathcal{C} \cap \mathcal{D}_{2^{-k} }} r(D)^s \leq \sum_{k=0}^K \sum_{D \in \mathcal{C}_j \cap \mathcal{D}_{2^{-k} }} r(D)^s \leq \varepsilon. \]
Letting $K \to \infty$ gives 
\[ \sum_{D \in \mathcal{C}} r(D)^s \leq \varepsilon. \]
So, $\cC$ satisfies condition $(2)$. By definition, it is easy to check $\cC_i\le \cC$ for every $\cC_i$ in the initial chain. This proves that $\cC$ is an upper bound of the chain. \end{proof}

\begin{remark}
{\rm
Besides $[0,1]^2$, this lemma holds for other compact metric spaces, for example $[0,1]^n$ or  $\S^2$. The proof is exactly the same.
}
\end{remark}

\begin{proof}[Proof that Theorem \ref{planediscrt} implies Theorem \ref{mainplane}]
Suppose $A\subset \R^3$ is a Borel set of Hausdorff dimension $\alpha$.
We may assume $A\subset B^3(0,1)$. Recall the definition of the exceptional set $$E_s:=\{\theta\in[0,1]:\dim \pi_\theta(A)< s\}.$$
Recall the definition of the $t$-dimensional Hausdorff content is given by
$$ \cH^t_\infty(B):=\inf\{ \sum_i r(B_i)^t: B\subset \cup_i B_i \}. $$
A property for the Hausdorff content is that 
$$ \dim (B)=\sup\{t:\cH^t_\infty(B)>0\}. $$
We choose $a<\dim(A), t<\dim(E_s)$. Then $\cH_\infty^t(E_s)>0$, and by Frostman's lemma there exists a probability measure $\nu_A$ supported on $A$ satisfying $\nu_A(B_r)\lesssim r^a$ for any $B_r$ being a ball of radius $r$. 
We may assume $t>0$, otherwise $\dim(E_s)=0$. We only need to prove
$$ a\le 1+s-t, $$
since then we can send $a\rightarrow \dim(A)$ and $t\rightarrow \dim(E_s)$. As $a$ and $t$ are fixed, we may assume $\cH_\infty^t(E_s)\sim 1$ is a constant.

Fix a $\theta\in E_s$. By definition, we have $\dim \pi_\theta(A)<s$. We also fix a small number $\en$ which we will later send to $0$.
By Lemma \ref{usefullemma}, we can find a covering of $\pi_\theta(A)$ by disks $\D_\theta=\{D\}$, each of which has radius $2^{-j}$ for some integer $j>|\log_2\en|$. We define $\mathbb D_{\theta,j}:=\{D\in\mathbb D_\theta: r(D)=2^{-j}\}$.
Lemma \ref{usefullemma} yields the following properties:
\begin{equation}\label{rsless1}
    \sum_{D\in\mathbb D_\theta}r(D)^s<1;
\end{equation}
For each $j$ and $r$-ball $B_r\subset V_\theta$ with $2^{-j}\le r\le 1$, we have
\begin{equation}\label{structure}
    \#\{D\in \D_{\theta,j}: D\subset B_r\}\lesssim (\frac{r}{2^{-j}})^s.
\end{equation}

For each $\theta\in E_s$, we can find such a $\D_\theta$. We also define the tube sets $\T_{\theta,j}:=\{\pi^{-1}_\theta(D): D\in\D_{\theta,j}\}\cap B^3(0,1)$, $\T_{\theta}=\bigcup_j\T_{\theta,j}$. Each tube in $\T_{\theta,j}$ has dimensions $2^{-j}\times 2^{-j}\times 1$ and direction $\gamma(\theta)$. One easily sees that $A\subset \bigcup_{T\in \T_\theta}T $. By pigeonholing, there exists $j(\theta)$ such that
\begin{equation}\label{pigeon1}
    \nu_A(A\cap(\cup_{T\in\T_{\theta,j(\theta)}}T ))\ge \frac{1}{10j(\theta)^2}\nu_A(A)=\frac{1}{10j(\theta)^2}.
\end{equation} 
For each $j>|\log_2\en|$, define $E_{s,j}:=\{\theta\in E_s: j(\theta)=j\}$. Then we obtain a partition of $E_s$:
$$ E_s=\bigsqcup_j E_{s,j}. $$
By pigeonholing again, there exists $j$ such that
\begin{equation}\label{pigeon2}
    \cH_\infty^t(E_{s,j})\ge \frac{1}{10j^2}\cH_\infty^t(E_s)\sim \frac{1}{10j^2}. 
\end{equation} 
In the rest of the poof, we fix this $j$. We also set $\de=2^{-j}(<\en)$. By Lemma \ref{lemdelta}, there exists a $\de$-separated $(\de,t)$-set $\Theta\subset E_{s,j}$ with cardinality $\#\Theta\gtrsim (\log\de^{-1})^{-2}\de^{-t}$.

Next, we consider the set $S:=\{(x,\theta)\in A\times \Theta: x\in\cup_{T\in\T_{\theta,j}}T \}$. We also use $\mu$ to denote the counting measure on $\Theta$.
Define the section of $S$:
$$ S_x=\{\theta: (x,\theta)\in S\},\ \ \  S_\theta:=\{x: (x,\theta)\in S\}. $$
By \eqref{pigeon1} and Fubini, we have
\begin{equation}\label{pigeon3}
    (\nu_A\times \mu)(S)\ge \frac{1}{10j^2} \mu(\Theta).
\end{equation}
This implies
\begin{equation}\label{pigeon4}
    (\nu_A\times \mu)\bigg(\Big\{(x,\theta)\in S: \mu(S_x)\ge\frac{1}{20j^2}\mu(\Theta)  \Big\}\bigg)\ge \frac{1}{20j^2} \mu(\Theta).
\end{equation} 
since
\begin{equation}
    (\nu_A\times \mu)\bigg(\Big\{(x,\theta)\in S: \mu(S_x)\le\frac{1}{20j^2}\mu(\Theta)  \Big\}\bigg)\le \frac{1}{20j^2}\cH_\infty^a(A) \mu(\Theta).
\end{equation} 
By \eqref{pigeon4}, we have
\begin{equation}\label{pigeon5}
    \nu_A\bigg(\Big\{x\in A: \mu(S_x)\ge \frac{1}{20j^2}\mu(\Theta) \Big\}\bigg)\ge \frac{1}{20j^2}. 
\end{equation} 

We are ready to apply Theorem \ref{planediscrt}. Recall $\de=2^{-j}$ and $\#\Theta\gtrsim (\log\de^{-1})^{-2}\de^{-t}$. By \eqref{pigeon5}, we can find a $\de$-separated subset of $\{x\in A: \# S_x\ge \frac{1}{20j^2}\#\Theta \}$ with cardinality $\gtrsim (\log\de^{-1})^{-2}\de^{-a}$. We denote the $\de$-neighborhood of this set by $H$, which is a union of $\de$-balls. For each $\de$-ball $B_\de$ contained in $H$, we see that there are $\gtrsim (\log\de^{-1})^{-2}\#\Theta$ many tubes from $\cup_{\theta\in\Theta}\T_{\theta,j}$ that intersect $B_\de$. We can now apply Theorem \ref{planediscrt} to obtain
$$ (\log\de^{-1})^{-4}\de^{-a-t}\lesssim\#\Theta\#H\le C_{s,\e}\de^{-1-s-\e}. $$
Letting $\en\rightarrow 0$ (and hence $\de\rightarrow 0$) and then $\e\rightarrow 0$, we obtain $a+t\le 1+s$.

\end{proof}

\subsection{Proof of Theorem \ref{planediscrt}}
The proof of Theorem \ref{planediscrt} is base on the $L^6$ decoupling inequality for cone which is well-understood.
For convenience, we will prove the following version of Theorem \ref{planediscrt} after rescaling $x\mapsto \de^{-1}x$.
\begin{theorem}\label{planerescale} Fix $0<s<2$. For each $\e>0$, there exists $C_{s,\e}$ so that the following holds. Let $\d>0$. Let $H\subset B^3(0,\de^{-1})$ be a union of $\de^{-a}$ many disjoint unit balls so that $H$ has measure $|H|\sim \de^{-a}$. Let $\Theta$ be a $\de$-separated subset of $[0,1]$ so that $\Theta$ is a $(\de,t)$-set and $\#\Theta\gtrsim (\log\de^{-1})^{-2} \de^{-t}$. Assume for each $\theta\in \Theta$, we have a collection of $1\times 1\times \de^{-1}$-tubes $\T_\theta$ pointing in direction $\ga(\theta)$. $\T_\theta$ satisfies the $s$-dimensional condition: 
\begin{enumerate}
    \item $\#\T_\theta\lesssim \de^{-s}$, 
    \item $\#\{T\in\T_\theta:T\cap B_r\}\lesssim r^{s}$, for any $B_r$ being a ball of radius $r$ $(1\le r\le \de^{-1})$.
\end{enumerate}
We also assume that each unit ball contained in $H$ intersects $\gtrsim |\log\de^{-1}|^{-2}\#\Theta$ many tubes from $\cup_\theta \T_\theta$.
Then 
\[ \d^{-t-a}\le C_{s,\e}\de^{-1-s-\e}.\] 
\end{theorem}

We first discuss the geometry of $\ga$.
Let $\ga$ be the non-degenerate curve as discussed in the beginning of this section. We have $|\ga'(\theta)|=1$. 
For convenience, we define 
\begin{equation}\label{coord}
    \be_1(\theta):=\ga(\theta),\ \  \be_2(\theta):=\ga'(\theta),\ \  \be_3(\theta):=\ga(\theta)\times \ga'(\theta).
\end{equation}
We see that $\{\be_1(\theta),\be_2(\theta),\be_3(\theta)\}$ form a Frenet coordinate along $\ga$.
Define the corresponding conical surface $\Ga:=\{r\be_3(\theta):1/2\le r\le 1, \theta\in[0,1]\}$. 

We first show that $\Ga$ satisfies the same non-degenerate condition as the standard cone. %, that is to say, $\Ga$ is generated by a non-degenerate planar curve. Equivalently, we just need to show that there is a plane $\Pi$ such that $\Ga\cap\Pi$ is a non-degenerate curve. %We write $\be_3(\theta)$ in coordinate as $(e_{31}(\theta),e_{32}(\theta),e_{33}(\theta))$, and by rotation we may assume $e_{33}(\theta)>0$ is strictly positive.
%Choosing $\Pi=\{x_3=1\}$, we see that $\Ga\cap\Pi$ is the curve $\wt\ga(\theta)=\frac{1}{e_{33}(\theta)}\be_3(\theta),\ \theta\in[0,1]$. To see $\wt\ga$ is non-degenerate, we just need to show the norm of $\wt\ga''(\theta)$ is strictly positive. 
Note that we have the following formulae for the Frenet coordinate:
\begin{align}
&\be_1'(\theta)=\be_2(\theta),\\
&\be_2'(\theta)=-\be_1(\theta)+\kappa(\theta)\be_3(\theta),\\
&\be_3'(\theta)=-\kappa(\theta)\be_2(\theta),
\end{align}
where $\kappa(\theta)=\langle\be_2'(\theta),\be_3(\theta)\rangle>0$. %After some computations, we have
%$$ \wt\ga''(\theta)=a(\theta)\be_2(\theta)+b(\theta)\be_3(\theta)+\frac{\kappa(\theta)}{e_{33}(\theta)}\be_1(\theta). $$
%Therefore, $|\wt\ga''(\theta)|\ge \frac{\kappa(\theta)}{e_{33}(\theta)}$ is strictly positive.

%Next, we
First, we show that $\Ga$ is a $C^2$ surface. We will do this by finding a reparametrization $s=s(\theta)$ so that $\be_3(\theta(s))$ is a $C^2$ function of $s$. Choose 
$$ s(\theta)=\int_0^\theta \kappa(t)dt, $$
and then $\frac{d\theta}{ds}=\kappa(\theta)^{-1}$. We have
$$ \frac{d}{ds}\be_3=\frac{d\theta}{ds}\cdot\frac{d}{d\theta}\be_3=-\be_2. $$
Since $\theta=\theta(s)$ is $C^1$, we have $\frac{d}{ds}\be_3=-\be_2(\theta(s))$ is $C^1$ with respect to $s$, and therefore $\be_3(\theta(s))$ is $C^2$ with respect to $s$. Moreover, $\det\left( \mathbf{e}_3(\theta(s)), \frac{d}{ds} \mathbf{e}_3(\theta(s)), \frac{d^2}{ds^2} \mathbf{e}_3(\theta(s)) \right) = \theta'(s)\det(\mathbf{e}_3(\theta(s)), -\mathbf{e}_2(\theta(s)), \mathbf{e}_1(\theta(s)))$ by the above, which is nonvanishing since $\theta'(s)$ is nonvanishing. 

\bigskip

For any large scale $R$, 
there is a standard partition of $N_{R^{-1}}\Ga$ into planks $\si$ of dimensions $R^{-1}\times R^{-1/2}\times 1$: $$N_{R^{-1}}\Ga=\bigcup \si.$$
For any Schwartz function $f$, we define $f_\si:=(1_\si\wh f)^\vee$ as usual. We have the following $L^6$-decoupling inequality for these planks. 

\begin{theorem}[Bourgain-Demeter \cite{bourgain2015proof}]\label{planedec}
For any Schwartz $f$ with $\wh f\subset N_{R^{-1}}\Ga$, we have
\begin{equation}\label{dec}
\|f\|_6\lesssim_\e R^{\e} \big(\sum_{\si:R^{-1}\times R^{-1/2}\times 1}\|f_\si\|_6^2\big)^{1/2}.    
\end{equation}
\end{theorem}

\begin{remark}\label{rmkcone}
We will actually apply Theorem \ref{planedec} to a slightly different cone
\begin{equation}
    \Ga_{K^{-1}}=\{r\be_3(\theta):K^{-1}\le r\le 1,\theta\in[0,1]\},
\end{equation}
for some $K\sim(\log\de^{-1})^{O(1)}$. Compared with $\Ga$, we see that $\Ga_{K^{-1}}$ is at distance $K^{-1}$ from the origin, but we still have a similar decoupling inequality.
Instead of \eqref{dec}, we have
\begin{equation}
    \|f\|_6\lesssim_\e K^{O(1)}R^\e \big(\sum_{\si:R^{-1}\times R^{-1/2}\times 1}\|f_\si\|_6^2\big)^{1/2}.
\end{equation}
The idea is to partition $\Ga_{K^{-1}}$ into $\sim O(K)$ many parts, each of which is roughly a cone for which we can apply Theorem \ref{planedec}. By triangle inequality, this results in an additional factor $K^{O(1)}$. It turns out that this factor is not harmful, since we will set $K\sim (\log R)^{O(1)}$ which can  be absorbed into $R^\e$. 

\end{remark}

We  are ready to prove Theorem \ref{planerescale}.
\begin{proof}[Proof of Theorem \ref{planerescale}] 
Recall that $\T_\theta$ is a collection of $1\times 1\times \de^{-1}$-tubes pointing to direction $\ga(\theta)=\be_1(\theta)$. 
We consider the dual of each $T_\theta$ in the frequency space.
For each $\theta$, we define $P_\theta$ to be a slab centered at the origin that has dimensions $1\times 1\times \de$, and its shortest direction is parallel to $\be_1(\theta)$. We see that $P_\theta$ is the dual rectangle of each $T_\theta\in\T_\theta$. Now, for each $T_\theta\in\T_\theta$, we choose a bump function $\psi_{T_\theta}$ satisfying the following properties: $\psi_{T_\theta}\ge 1$ on $T_\theta$, $\psi_{T_\theta}$ decays rapidly outside $T_\theta$, and $\supp \wh\psi_{T_\theta}\subset P_\theta$.

Define functions 
\[ f_\theta=\sum_{T_\theta\in\T_\theta}\s_{T_\theta}\qquad\text{and}\qquad f=\sum_{\theta\in \Theta}f_\theta. \]
From our definitions, we see that for any $x\in H$, we have $f(x)\gtrsim \#\{T\in\cup_\theta \T_\theta: x\in 
T\} \gtrsim (\log\de^{-1})^{-2}\#\Theta$. Therefore, we obtain
\begin{equation}
    |H| (\log\de^{-1})^{-2p}(\#\Theta)^p \le 2^p\int_{H}|f|^p,
\end{equation}
for any $p$. For our purpose, we just choose $p=6$, so we have
\begin{equation}\label{plane1}
    |H| (\log\de^{-1})^{-12}(\#\Theta)^6\le 2^6\int_{H}|f|^6,
\end{equation}
Our goal is to find an upper bound for the right hand side of \eqref{plane1}. We will decompose $P_\theta$ into pieces and estimate the contribution of $\wh f_\theta$ from each piece.

Let us discuss the decomposition for $P_\theta$. Recall that $P_\theta$ is a $1\times 1\times \de$-slab centered at the origin with normal direction $\be_1(\theta)$.  
Recall \eqref{coord}, we can write $P_\theta=\{\sum_{i=1}^3\xi_i\be_i(\theta): |\xi_1|\le \de,|\xi_2|\le 1, |\xi_3|\le 1\}$.

%%%%%%%%%%%%%%%%%
\begin{figure}%theta-omega
\begin{tikzpicture}

\begin{scope}
   [x={(6cm,0)},
    y={({cos(45)*.5cm},{sin(45)*.5cm})},
    z={({cos(90)*6cm},{sin(90)*6cm})},line join=round, black]
  \draw[] (0,0,0) -- (0,0,1) -- (0,1,1);
   \draw[] (0,0,0) -- (1,0,0) -- (1,1,0);

  \draw[] (1,0,0) -- (1,0,1) -- (1,1,1) -- (1,1,0) -- cycle;
 \draw[] (0,0,1) -- (1,0,1) -- (1,1,1) -- (0,1,1) -- cycle;
 \end{scope}  
 
%low
\begin{scope}
   [shift={(1.5cm,1.5cm)},x={(3cm,0)},
    y={(0,0)},
    z={({cos(90)*3cm},{sin(90)*3cm})},line join=round, black]
  \draw[] (0,0,0) -- (0,0,1) -- (1,0,1) -- (1,0,0) --cycle;
\end{scope}

%left high
\begin{scope}
   [x={(1.5cm,0)},
    y={({cos(45)*.5cm},{sin(45)*.5cm})},
    z={({cos(90)*6cm},{sin(90)*6cm})},line join=round, black]
  \draw[] (0,0,0) -- (0,0,1) -- (0,1,1);
   \draw[] (0,0,0) -- (1,0,0);
  \draw[] (1,0,0) -- (1,0,1);
 \draw[] (0,0,1) -- (1,0,1) -- (1,1,1) -- (0,1,1) -- cycle;
 \end{scope}  

%right high
\begin{scope}
   [shift={(4.5cm,0)},x={(1.5cm,0)},
    y={({cos(45)*.5cm},{sin(45)*.5cm})},
    z={({cos(90)*6cm},{sin(90)*6cm})},line join=round, black]
  \draw[] (0,0,0) -- (0,0,1) -- (0,1,1);
   \draw[] (0,0,0) -- (1,0,0);
  \draw[] (1,0,0) -- (1,0,1);
 \draw[] (0,0,1) -- (1,0,1) -- (1,1,1) -- (0,1,1) -- cycle;
 \end{scope}

%\lambda=\delta^{1/2}
\begin{scope}
   [shift={(2.75cm,0)},x={(0.5cm,0)},
    y={({cos(45)*.5cm},{sin(45)*.5cm})},
    z={({cos(90)*6cm},{sin(90)*6cm})},line join=round, black]
  \draw[] (0,0,0.75) -- (0,0,1) -- (0,1,1);
   \draw[] (0,0,0) -- (0,0,0.25);
  \draw[] (1,0,0) -- (1,0,0.25);
  \draw[] (1,0,0.75) -- (1,0,1);
 \draw[] (0,0,1) -- (1,0,1) -- (1,1,1) -- (0,1,1) -- cycle;
 
 \end{scope} 

\node[scale=0.8] at (0.8,3)
{$P_{\theta,high}$};

\node[scale=0.8] at (3,3)
{$P_{\theta,low}$};

\node[scale=0.8] at (3,5.2)
{$P_{\theta,\delta^{1/2}}$};

\node[scale=0.8] at (2,5.2)
{$P_{\theta,\lambda}$};

\draw[black, ->] (10,0) -- (11,0) node[right]{$\boldsymbol e_2(\theta)$};
\draw[black, ->] (10,0) -- (9.5,-0.8) node[below]{$\boldsymbol e_1(\theta)$};
\draw[black, ->] (10,0) -- (10,1) node[above]{$\boldsymbol e_3(\theta)$};

% \node[scale=0.8] at (1.38,4.7) {$\theta_{\tau_2,high}^+$};
% \node[scale=0.8] at (.8,5.2) {$\theta_{\tau_1,high}^+$};
% \node at (.8,6) {$\omega$}; 
 
\end{tikzpicture}
\caption{High-low decomposition for $P_\theta$}
\label{decomposition}
\end{figure}
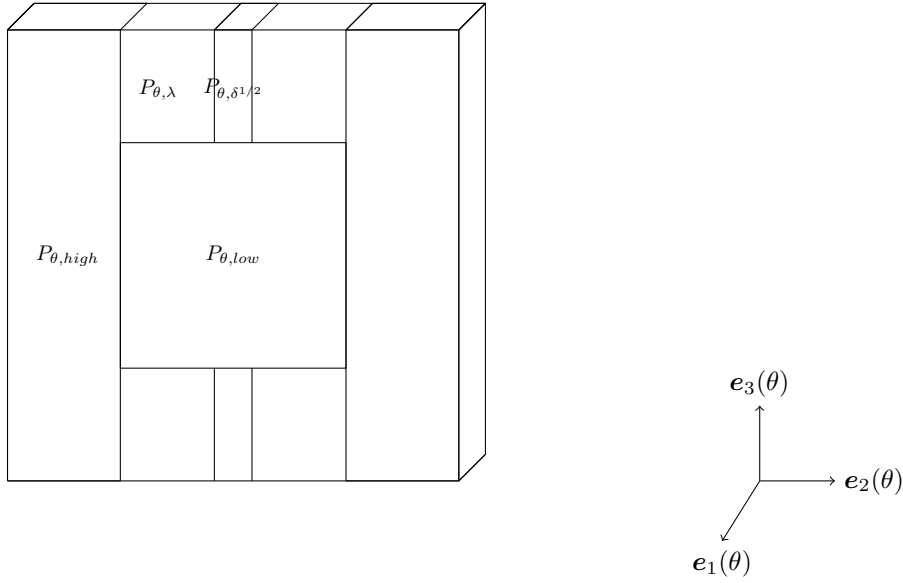

\begin{definition}\label{defdecomposition}
See Figure \ref{decomposition}.
Let $K$ be a large number which we will choose later. (Actually, we will choose $K\sim (\log\de^{-1})^{O(1)}$) Define the high part of $P_\theta$ as
$$ P_{\theta,high}:=\{\sum_{i=1}^3\xi_i\be_i(\theta): |\xi_1|\le \de, K^{-1}\le|\xi_2|\le 1, |\xi_3|\le 1\}. $$
Define the low part of $P_\theta$ as 
$$ P_{\theta,low}:=\{\sum_{i=1}^3\xi_i\be_i(\theta): |\xi_1|\le \de, |\xi_2|\le K^{-1}, |\xi_3|\le K^{-1}\}. $$
For dyadic numbers $\lambda\in(\de^{1/2},K^{-1}]$, define 
$$ P_{\theta,\lambda}:=\{\sum_{i=1}^3\xi_i\be_i(\theta): |\xi_1|\le \de, \frac{1}{2}\lambda\le|\xi_2|\le \lambda, K^{-1}\le|\xi_3|\le 1\}. $$
In particular, we define
$$ P_{\theta,\de^{1/2}}:=\{\sum_{i=1}^3\xi_i\be_i(\theta): |\xi_1|\le\de, |\xi_2|\le \de^{1/2},K^{-1}\le|\xi_3|\le 1  \}. $$
\end{definition}

\begin{remark}
We obtain a partition of $P_\theta$ as 
$$P_\theta=P_{\theta,high}\bigsqcup P_{\theta,low}\bigsqcup_\lambda P_{\theta,\lambda}.$$
We see that $P_{\theta,\lambda}$ consists of four planks of dimensions  $\sim \de\times\lambda\times 1$ whose longest side is along direction $\be_3(\theta)$. Here $\lambda$ plays a role of angular parameter in the sense that  $P_{\theta,\lambda}$ are roughly those points in $P_{\theta}\setminus P_{\theta,low}$ so that the lines connecting them with the origin form an angle $\sim \lambda $ with $\be_3(\theta)$.
\end{remark}

We choose a smooth partition of unity adapted to this covering which we denote by $\eta_{\theta,high}, \eta_{\theta,low}, \eta_{\theta,\lambda}$, so that
\begin{equation}
    \eta_{\theta,high}+\eta_{\theta,low}+ \sum_{\de^{1/2}\le\lambda\le K^{-1}}\eta_{\theta,\lambda}=1
\end{equation}
on $P_\theta$. Since $\supp\wh f_\theta\subset P_\theta$, we also obtain a decomposition of $f_\theta$
\begin{equation}\label{decompftheta}
    f_\theta=f_{\theta,high}+f_{\theta,low}+ \sum_{\de^{1/2}\le\lambda\le K^{-1}}f_{\theta,\lambda},
\end{equation}
where $\wh f_{\theta,high}=\eta_{\theta,high} \wh f_\theta, \wh f_{\theta,low}=\eta_{\theta,low}\wh f_{\theta}, \wh f_{\theta,\lambda}=\eta_{\theta,\lambda}\wh f_{\theta}.$
Similarly, we have a decomposition of $f$
\begin{equation}\label{decompf}
    f=f_{high}+f_{low}+ \sum_{\de^{1/2}\le\lambda\le K^{-1}}f_{\lambda},
\end{equation}
where $f_{high}=\sum_{\theta}f_{\theta,high},  f_{low}=\sum_\theta f_{\theta,low}, f_{\lambda}=\sum_\theta f_{\theta,\lambda}.$

Recalling \eqref{plane1} and using triangle inequality, we have
\begin{equation}\label{3cases}
    |H|(\log\de^{-1})^{-12}(\#\Theta)^6\lesssim \int_H|f_{low}|^6+ \int_H|f_{high}|^6+(\log\de^{-1})^{O(1)}\sum_{\lambda}\int_H|f_\lambda|^6. 
\end{equation}
We will discuss three cases depending on which term on the right hand side of \eqref{3cases} dominates.

\medskip

\noindent\fbox{Case 1: Low case} 

If the first term on the right hand side of \eqref{3cases} dominates, we say we are in the low case.
Actually, we will see that we are never in the low case by showing 
\begin{equation}\label{lowsmall}
    \int_H |f_{low}|^6\le C^{-1} |H|(\log\de^{-1})^{-12}(\#\Theta)^6,
\end{equation}
for some large constant $C$.
This means the low term on the right hand side of \eqref{3cases} will not dominate.
By properly choosing $K$, we can show a pointwise bound for $f_{low}$:
\begin{equation}\label{pointwise}
    |f_{low}(x)|\le C^{-1} (\log\de^{-1})^{-2}\#\Theta.
\end{equation}
This will immediately imply \eqref{lowsmall}. Let us focus on \eqref{pointwise}.

Recall that $f_{low}=\sum_{\theta}f_{\theta,low}=\sum_{\theta}f_{\theta}*\eta_{\theta,low}^\vee$.
Since $\eta_{\theta,low}$ is a bump function at $P_{\theta,low}$, we see that $\eta^\vee_{\theta,low}$ is a bump function essentially supported in the dual of $P_{\theta,low}$. Denote the dual of $P_{\theta,low}$ by $T_{\theta,K}$ which is a $K\times K\times \de^{-1}$-tube parallel to $\be_1(\theta)$. One has
$$ |\eta_{\theta,low}^\vee|\lesssim \frac{1}{|T_{\theta,K}|}\psi_{T_{\theta,K}}. $$
Here $\psi_{T_{\theta,K}}$ is a bump function $=1 $ on $T_{\theta,K}$ and decays rapidly outside $T_{\theta,K}$.

By definition, $f_\theta=\sum_{T_\theta}\psi_{T_\theta}$ is the sum of bump function of tubes. We have
\begin{align}
    |f_{low}|\lesssim\sum_{\theta}\sum_{T_\theta\in\T_\theta}\psi_{T_\theta}*\frac{1}{|T_{\theta,K}|}\psi_{T_{\theta,K}}.
\end{align}
If we ignore the rapidly decaying tails, we have
\begin{equation}
    |f_{low}(x)|\lesssim \sum_\theta\frac{1}{K^2} \#\{T_\theta\in\T_\theta:T_\theta\cap B_{100K}(x)\not=\emptyset\}.
\end{equation}
Recalling the condition (2) in Theorem \ref{planerescale}, we have $$\#\{T_\theta\in\T_\theta:T_\theta\cap B_{100K}(x)\not=\emptyset\}\lesssim (100K)^s.$$
This implies
\begin{equation}
    |f_{low}(x)|\lesssim  K^{s-2}\#\Theta.
\end{equation}
Since $s<2$, by choosing $K\sim (\log\de^{-1})^{\frac{2}{2-s}}$, we obtain \eqref{pointwise}.

In the rest of the proof, we may pretend $K$ is a large constant, since any $(\log\de^{-1})^{O(1)}$-loss is allowable (see Remark \ref{rmkcone}).
\medskip

\noindent\fbox{Case 2: High case}

If the second term on the right hand side of \eqref{3cases} dominates, we say we are in the high case.

Since for any $x\in\R^3$ there is at most $1$ tube in $\T_\theta$ pass through $x$, we have $|f_\theta(x)|\lesssim 1$ (here we use $\lesssim 1$ instead of $\le 1$ to take care of the rapidly decaying tail). Recalling the definition $f_{high}=\sum_\theta \eta^\vee_{\theta,high}*f_\theta$ and noting that each $\eta^\vee_{\theta,high}$ is $L^1$ bounded, we have
$$ |f_{high}(x)|\lesssim \sum_{\theta}|f_\theta(x)|\le \#\Theta. $$
We see that 
\begin{equation}
    \int_H|f_{high}|^6\lesssim (\#\Theta)^4\int|f_{high}|^2=(\#\Theta)^4\int|\sum_\theta f_{\theta,high}|^2.
\end{equation}
Next we will show that $\{\supp \wh f_{\theta,high}\}_{\theta\in\Theta}$ are finitely overlapping, i.e., $\{P_{\theta,high}\}_{\theta\in\Theta}$ are finitely overlapping. (Actually they are $O(K)$-overlapping. But since the $K^{O(1)}$-loss are acceptable, we may just pretend $K\lesssim 1$. See also Remark \ref{rmkcone}.) If this is true, then we have
\begin{equation}
    \int_H|f_{high}|^6\lesssim (\#\Theta)^4\int \sum_\theta|f_{\theta,high}|^2.
\end{equation}
Since $\int \sum_\theta|f_{\theta,high}|^2=\int \sum_\theta|\eta^\vee_{\theta,high}*f_{\theta}|^2\le\int \sum_\theta |f_\theta|^2\sim \sum_\theta(\#\T_\theta) \de^{-1}.$ We obtain
\begin{equation}\label{highpart}
    |H|(\log\de^{-1})^{-12}(\#\Theta)^6\lesssim \int|f_{high}|^6\lesssim (\#\Theta)^4\sum_\theta(\#\T_\theta) \de^{-1}\lesssim (\#\Theta)^4 \de^{-s-t-1},
\end{equation}
which implies 
\begin{equation}
    \de^{-a-t}\lesssim \de^{-s-1-\e}.
\end{equation}
Now we prove that $\{P_{\theta,high}\}_{\theta\in\Theta}$ are finitely overlapping. 
First, recall that $$P_{\theta,high}=\{\sum_{i=1}^3\xi_i\be_i(\theta): |\xi_1|\le \de, K^{-1}\le|\xi_2|\le 1, |\xi_3|\le 1\}.$$
We see that $P_{\theta,high}$ is contained in the $\de$-neighborhood of the plane
$$ \Pi_{\theta,high}=\{\xi_2 \be_2(\theta)+\xi_3 \be_3(\theta): K^{-1}\le|\xi_2|\le 1, |\xi_3|\le 1\}. $$
To show the finitely overlapping property, we just need to show: For any $\theta\in[0,1]$ and $\theta'=\theta+\Delta\in [0,1]$ with $C\de\le\Delta\le C^{-1}$ (for some bounded $C$ to be determined later), if $\xi\in \Pi_{\theta,high}$, then 
$$ \dist(\xi,\Pi_{\theta+\Delta,high})>10\de. $$
Write $\xi=a\be_2(\theta)+b\be_3(\theta)= a\ga'(\theta)+b\ga(\theta)\times\ga'(\theta)$, where $|a|\in[K^{-1},1]$ and $|b|\le 1$. Since the normal direction of $\Pi_{\theta+\Delta,high}$ is $\ga(\theta+\Delta)$, it suffices to prove
\begin{equation}\label{disjoint}
    |\ga(\theta+\Delta)\cdot \big( a\ga'(\theta)+b\ga(\theta)\times\ga'(\theta) \big)|\ge 10\de.
\end{equation}
By Taylor's expansion, we have
$\ga(\theta+\Delta)=\ga(\theta)+\Delta\ga'(\theta)+O(\Delta^2)$.
We see the left hand side of \eqref{disjoint} is $ |a\Delta|\ga'(\theta)|^2+O(\Delta^2)|\ge |\big(a-O(\Delta)\big)\Delta|\ge 10\de $, if $C$ is large enough (depending on $K$).

\medskip

\noindent\fbox{Case 3: $\lambda$-middle case $(\d^{1/2}\le\lambda\le K^{-1})$}
If the term $(\log\de^{-1})^{O(1)}\sum_\lambda\int_H|f_\lambda|^6$ on the right hand side of \eqref{3cases} dominates, we say we are in the $\lambda$-middle case. We remark that when $\lambda$ is close to $K^{-1}$, $\int_H|f_\lambda|^6$ can be estimated in a similar way as in the High case. 
We will be interested in the cone
$$\Ga_{K^{-1}}=\{r\be_3(\theta):K^{-1}\le r\le 1, \theta\in[0,1]\}.$$
Recall Remark \ref{rmkcone} that we still have the decoupling inequality for this cone.

We first discuss the case that $\lambda=\de^{1/2}$.

\medskip

\noindent\fbox{Case 3.1: $\lambda =\d^{1/2}$} 

When $\lambda=\de^{1/2}$, we have $f_{\de^{1/2}}=\sum_{\theta}f_{\theta,\de^{1/2}}$,
where each $\wh f_{\theta,\de^{1/2}}$ is supported in $ P_{\theta,\de^{1/2}}$. Note that $P_{\theta,\de^{1/2}}$ consists of two pieces: One is $$ P_{\theta,\de^{1/2}}^+:=\{\sum_{i=1}^3\xi_i\be_i(\theta): |\xi_1|\le\de, |\xi_2|\le \de^{1/2},K^{-1}\le\xi_3\le 1  \}. $$
the other is 
$$ P_{\theta,\de^{1/2}}^-:=\{\sum_{i=1}^3\xi_i\be_i(\theta): |\xi_1|\le\de, |\xi_2|\le \de^{1/2},-1\le\xi_3\le -K^{-1}  \}. $$
We note that $P_{\theta,\de^{1/2}}^+$ lies in the $\de$-neighborhood of $\Ga_{K^{-1}}$. Symmetrically, $P_{\theta,\de^{1/2}}^-$ lies in the $\de$-neighborhood of $\Ga_{K^{-1}}^-$, which is the reflection of $\Ga_{K^{-1}}$ with respect to the origin.
We can write $f_{\theta,\de^{1/2}}=f_{\theta,\de^{1/2}}^++f_{\theta,\de^{1/2}}^-$, so that $\supp \wh f_{\theta,\de^{1/2}}^+\subset P_{\theta,\de^{1/2}}^+$ and $\supp \wh f_{\theta,\de^{1/2}}^-\subset P_{\theta,\de^{1/2}}^-$. We also write $f_{\de^{1/2}}^+=\sum_\theta f_{\theta,\de^{1/2}}^+$ and $ f_{\de^{1/2}}^-=\sum_\theta f_{\theta,\de^{1/2}}^-$, and then $f_{\de^{1/2}}=f_{\de^{1/2}}^++f_{\de^{1/2}}^-$.
We have
\begin{equation}
    \int |f_{\de^{1/2}}|^6\lesssim \int |f_{\de^{1/2}}^+|^6+\int |f_{\de^{1/2}}^-|^6.
\end{equation}
By symmetry, we only estimate $\int |f_{\de^{1/2}}^+|^6$.

Note that $P_{\theta,\de^{1/2}}^+$ and $P_{\theta',\de^{1/2}}^+$ are essentially the same when $|\theta-\theta'|\lesssim \de^{1/2}$; $P_{\theta,\de^{1/2}}^+$ and $P_{\theta',\de^{1/2}}^+$ are essentially distinct when $|\theta-\theta'|\gtrsim \de^{1/2}$. We can choose a partition of $N_\de(\Ga_{K^{-1}})$ by finitely overlapping planks of dimensions $10\de\times 10\de^{1/2}\times 1$, denoted by $\{R\}$. We attach each $P_{\theta,\de^{1/2}}^+$ to one of the $R$ and denote by $\theta\prec R$, if $P_{\theta,\de^{1/2}}^+\subset R$. We see that for each $R$, there are $\lesssim \de^{-1/2}$ many $P_{\theta,\de^{1/2}}$ attached to it. We define $f_R:=\sum_{\theta\prec R}f_{\theta,\de^{1/2}}^+$.
The Fourier support of $f_R$ is contained in $R$, by Theorem \ref{planedec}, we have
$$ \int |f_{\de^{1/2}}^+|^6=\int |\sum_R f_R|^6\lesssim_\e \de^{-\e}  (\#\{R: f_R\neq 0\})^2\sum_R  \|f_R\|^6_6 . $$
By pigeonholing, we may pass to a subset of $\{R\}$ so that $\#\{\theta:\theta\prec R\}$ are all comparable. For simplicity, we write $\#\{\theta:\theta\prec R\}$ as $\#\{\theta\prec R\}$, and write the number of $\{R\}$ after pigeonholing as $\#\{R\}$.
For each $R$, we have by triangle inequality:
$$ \|f_R\|_6^6\lesssim \#\{\theta\prec R\}^5\sum_{\theta\prec R}\|f_{\theta,\de^{1/2}}^+\|_6^6\lesssim \#\{\theta\prec R\}^5\sum_{\theta\prec R}\|f_{\theta}\|_6^6.$$
We obtain that
$$ \int|f_{\delta^{1/2}}^+|^6\lesssim_\e \de^{-\e} \#\{R\}^2\#\{\theta\prec R\}^5\sum_\theta\|f_\theta\|_6^6. $$
Note that $\#\{R\}\#\{\theta\prec R\}\lesssim \#\Theta$, $\#\{\theta\prec R\}\lesssim \de^{-t/2}$ (by the $(\de,t)$-spacing of $\Theta$), and $\sum_\theta\|f_\theta\|_6^6\lesssim \#\Theta \de^{-1-s}$. We obtain
$$ \int|f_{\delta^{1/2}}^+|^6\lesssim_\e \de^{-\e} (\#\Theta)^3\de^{-3t/2-1-s}. $$
Similarly, we have
$$ \int|f_{\delta^{1/2}}^-|^6\lesssim_\e \de^{-\e} (\#\Theta)^3\de^{-3t/2-1-s}. $$
As a result, we obtain
\begin{equation}
    |H|(\#\Theta)^6\lesssim_\e \de^{-2\e} (\#\Theta)^3\de^{-3t/2-1-s}.
\end{equation}
Combined with $\#\Theta\gtrsim (\log\de^{-2})^{-2}\de^{-t}$, we have $\de^{-a-3t/2}\lesssim \de^{-\e-1-s}$ which is even better than what we aimed.

\medskip

\noindent\fbox{Case 3.2: $\lambda \in (\de^{1/2},K^{-1}]$}

For $\lambda$ being a dyadic scale in $(\de^{1/2},K^{-1}]$, we see that $\wh f_{\theta,\lambda}$ is supported in $P_{\theta,\lambda}$ which consists of four separated planks ($P_{\theta,\de^{1/2}}$ only consists of two planks). As in the proof of case $\lambda=\de^{1/2}$, we will write $f_{\theta,\lambda}$ as the sum of four functions each of which has Fourier support in one of the planks of $P_{\theta,\lambda}$. We will estimate for one of the planks. 
We define
\begin{equation}
    Q_{\theta}:=\{\sum_{i=1}^3\xi_i\be_i(\theta): |\xi_1|\le\de, \frac{1}{2}\lambda\le\xi_2\le \lambda, K^{-1}\le \xi_3\le 1\}.
\end{equation}
% For convenience, we also define
% $$ P_{\theta}:=\{\xi_1\ga'(\theta)+\xi_2\ga(\theta)\times\ga'(\theta)+\xi_3\ga(\theta): |\xi_1|\le \lambda\de^{-1}, K^{-1}\de^{-1}\le \xi_2\le \de^{-1},|\xi_3|\le 1  \}. $$
Roughly speaking, $Q_{\theta}$ is the top-right plank of $P_{\theta,\lambda}$ and the distance between $Q_{\theta}$ and the line $\R\be_3(\theta)$ is $\ge \frac{1}{2}\lambda$. For simplicity, we may assume $f_{\theta,\lambda}$ has Fourier support in $Q_\theta$.

We discuss some geometric properties for the planks $\{Q_{\theta}\}_{\theta\in\Theta}$. First of all, there is a canonical finitely overlapping covering of $N_{\lambda^2}(\Ga_{K^{-1}})$ by planks of dimensions $\lambda^2\times \lambda\times 1$. 
More precisely, we choose $\Si=\lambda\ZZ\cap[0,1]$ to be a set of $\lambda$-lattice points. For each $\si\in\Si$, define
$$ R_\si:=\{\sum_{i=1}^3\xi_i\be_i(\si): |\xi_1|\le C_1\lambda^2, |\xi_2|\le C_1\lambda, C_1^{-1}K^{-1}\le \xi_3\le C_1  \}, $$
where $C_1$ is a large constant.
We see that $\{R_\si\}$ form a finitely overlapping covering of $N_{\lambda^2}(\Ga_{K^{-1}})$. We have the following three properties:

\begin{lemma} \label{disjointlemma}
For $Q_\theta, R_\si$ defined above, we have
\begin{enumerate}
\item If $|\theta-\si|\lesssim \lambda$, then $Q_{\theta}$ is contained in $R_\si$.\label{property1}
\item If $|\theta-\theta'|\lesssim \lambda^{-1}\de$, then $Q_{\theta}$ and $Q_{\theta'}$ are essentially the same. \label{property2}
\item If $|\theta-\theta'|\gtrsim \lambda^{-1}\de$, then $Q_{\theta}$ and $Q_{\theta'}$ are disjoint.\label{property3}
\end{enumerate}
\end{lemma}

Before proving the lemma, we see how it can be used to finish the proof of Theorem \ref{planerescale}.
Motivated by Property \eqref{property2} and \eqref{property3}, we define $\Tau=(\lambda^{-1}\de)\ZZ\cap[0,1]$, and for each $\tau\in\Tau$ define
$$ S_\tau:=\{\sum_{i=1}^3 \xi_i\be_i(\tau): |\xi_1|\le \de, |\xi_2|\le C_2\lambda, C_2^{-1}K^{-1}\le \xi_3\le C_2\}, $$
where $C_2$ is a large constant but much smaller than $C_1$. Note that $S_\tau$ has the same dimensions as $Q_\theta$ up to a $C_2$-dilation.

Now we have three subsets of $[0,1]$: $$\Theta=\de\ZZ\cap[0,1],\ \  \Tau=(\lambda^{-1}\de)\ZZ\cap[0,1],\ \  \Si=\lambda\ZZ\cap[0,1].$$
We will define a relationship between their elements. For any $\theta\in\Theta$, we attach it to a $\tau\in\Tau$ such that $|\theta-\tau|\le \lambda^{-1}\de$, which we denote by $\theta\prec\tau$. For any $\tau\in\Tau$, we attach it to a $\si\in\Si$ such that $|\tau-\si|\le \lambda$, which we denote by $\tau\prec \si$. We also write $\theta\prec\si$ if there is a $\tau$ such that $\theta\prec\tau$ and $\tau\prec \si$.

By property \eqref{property1}, if $\theta\prec\si$, then $Q_\theta\subset R_\si$. By property \eqref{property2}, for a given $\tau\in\Tau$, all the planks in $\{Q_{\theta}:\theta\prec \tau\}$ are essentially the same and contained in $S_\tau$. By property \eqref{property3}, if $Q_\theta, Q_{\theta'}$ lie in different $S_\tau$, then $Q_\theta, Q_{\theta'}$ are disjoint. 

Before estimating $\int |f_\lambda|^6$, we may apply a pigeonhole argument to pass to subsets of $\Theta,\Tau,\Sigma$ (still denoted by $\Theta,\Tau,\Sigma$), so that $\#\{\theta\in\Theta:\theta\prec \tau\}$ are comparable for $\tau\in\Tau$, and $\#\{\tau\in\Tau:\tau\prec \si  \}$ are comparable for $\si\in\Sigma$. For convinience, we write $\#\{\theta\in\Theta:\theta\prec \tau\}$ as $\#\{\theta\prec \tau\}$, and write $\#\{\tau\in\Tau:\tau\prec \si  \}$ as $\#\{\tau\prec \si  \}$.

We define 
$$f_\tau:=\sum_{\theta\prec\tau}f_{\theta,\lambda},$$
$$ f_\si:=\sum_{\theta\prec\si}f_{\theta,\lambda}=\sum_{\tau\prec\si}f_\tau. $$
By decoupling for $N_{\lambda^2}(\Ga_{K^{-1}})=\bigsqcup_\si R_\si$, we have
\begin{equation}
    \int |f_\lambda|^6=\int |\sum_\si f_\si|^6\lesssim \de^{-\e}\#\{\si\}^2 \sum_\si \int |f_\si|^6.
\end{equation}

By the trivial decoupling for $R_\si=\bigsqcup_{\tau\prec\si}S_\tau$, we have
\begin{equation}
    \int |f_\si|^6=\int |\sum_{\tau\prec\si}f_\tau|^6\lesssim \#\{\tau\prec\si\}^{4}\sum_{\tau\prec\si}\int |f_\tau|^6.
\end{equation}

By H\"older's inequality, we have
\begin{equation}
    \int |f_\tau|^6=\int |\sum_{\theta\prec\tau}f_{\theta,\lambda}|^6\lesssim \#\{\theta\prec\tau\}^5\sum_{\theta\prec\tau}\int |f_{\theta,\lambda}|^6.
\end{equation}

Combining the three inequalities, we obtain
\begin{equation}
    \int|f_\lambda|^6\lesssim \de^{-\e}\#\{\si\}^2\#\{\tau\prec\si\}^4\#\{\theta\prec\tau\}^5\sum_\theta\int|f_{\theta,\lambda}|^6.
\end{equation}
Note that $\#\{\si\}\#\{\tau\prec \si\}\#\{\theta\prec\tau\}\lesssim \#\Theta$, $\#\{\tau\prec \si\}\#\{\theta\prec\tau\}\lesssim (\lambda\de^{-1})^t$ (by the $(\de,t)$-spacing of $\Theta$), and $\#\{\theta\prec\tau\}\lesssim \lambda^{-t}$ (by the $(\de,t)$-spacing of $\Theta$). We also note that $\sum_\theta\int|f_{\theta,\lambda}|^6\lesssim (\#\T) \de^{-1} \lesssim \#\Theta \de^{-s-1}$.
We obtain
\begin{equation}\label{lambdapart}
    \int|f_\lambda|^6\lesssim \de^{-\e}(\#\Theta)^2(\lambda\de^{-1})^{2t}\lambda^{-t}(\#\T)\de^{-1}\lesssim\de^{-\e}(\#\Theta)^3\lambda^t\de^{-s-1-2t}.
\end{equation}

Plugging into \eqref{3cases} and noting $\de^{-t}\gtrsim\#\Theta \gtrsim (\log\de^{-1})^{-2}\de^{-t}$, we obtain
$$ \de^{-a-t}\lesssim \de^{-2\e} \lambda^t \de^{-s-1}, $$
which is better than we aimed because of the factor $\lambda^t$.

\end{proof}

%%%%%%%%%%%%%%%%%%%%%
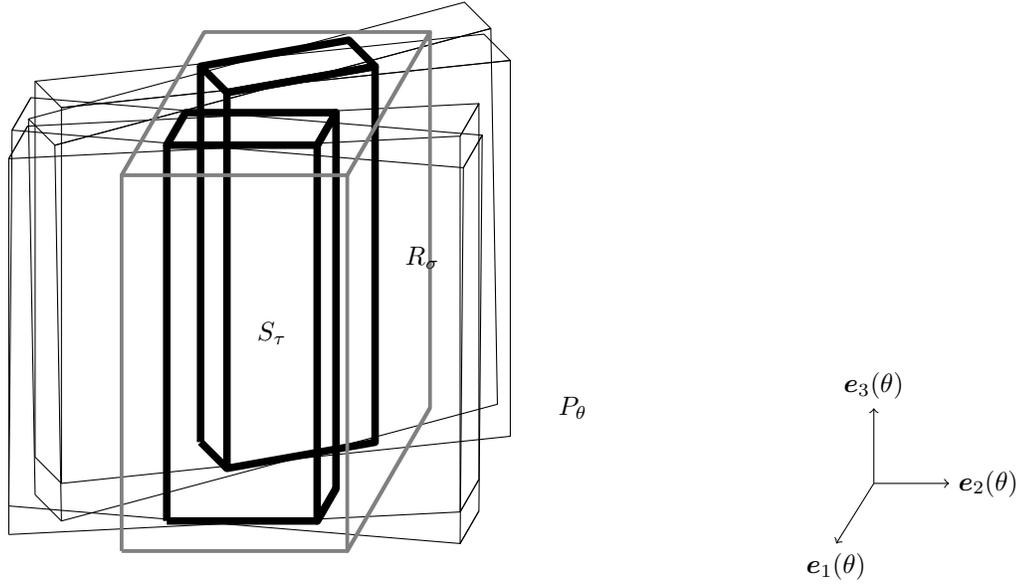
\begin{figure}%theta-omega
\begin{tikzpicture}
%first small plank(bottom)
 \begin{scope}
   [shift={(0.6,-.5)},x={(2cm,0cm)},
    y={({cos(60)*.5cm},{sin(60)*.5cm})},
    z={({cos(90)*5cm},{sin(90)*5cm})},line join=round,black, line width=1mm]
  \draw[] (0,0,0) -- (0,0,1) -- (0,1,1);
   \draw[] (0,0,0) -- (1,0,0) -- (1,1,0);

  \draw[] (1,0,0) -- (1,0,1) -- (1,1,1) -- (1,1,0) -- cycle;
 \draw[] (0,0,1) -- (1,0,1) -- (1,1,1) -- (0,1,1) -- cycle;
 \end{scope}

% normal caps
\begin{scope}
   [shift={(-1.5cm,-.68)},x={(6cm,0.3cm)},
    y={({cos(60)*.5cm},{sin(60)*.5cm})},
    z={({cos(90)*5cm},{sin(90)*5cm})},line join=round, black]
  \draw[] (0,0,0) -- (0,0,1) -- (0,1,1);
   \draw[] (0,0,0) -- (1,0,0) -- (1,1,0);

  \draw[] (1,0,0) -- (1,0,1) -- (1,1,1) -- (1,1,0) -- cycle;
 \draw[] (0,0,1) -- (1,0,1) -- (1,1,1) -- (0,1,1) -- cycle;
 \end{scope} 

\begin{scope}
   [shift={(-1.5cm,-.3)},x={(6cm,-0.5cm)},
    y={({cos(60)*.5cm},{sin(60)*.5cm})},
    z={({cos(89.5)*5cm},{sin(89.5)*5cm})},line join=round, black]
  \draw[] (0,0,0) -- (0,0,1) -- (0,1,1);
   \draw[] (0,0,0) -- (1,0,0) -- (1,1,0);

  \draw[] (1,0,0) -- (1,0,1) -- (1,1,1) -- (1,1,0) -- cycle;
 \draw[] (0,0,1) -- (1,0,1) -- (1,1,1) -- (0,1,1) -- cycle;
 \end{scope}

 %second small plank
\begin{scope}
   [shift={(1.4,.2)},x={(-.35cm,.35cm)},
    y={({cos(10)*2cm},{sin(10)*2cm})},
    z={({cos(90)*5cm},{sin(90)*5cm})},line join=round,fill opacity=0.5,black, line width=1mm]
  \draw[] (0,0,0) -- (0,0,1) -- (0,1,1) -- (0,1,0) -- cycle;
   \draw[] (0,0,0) -- (1,0,0);
   \draw[] (1,0,0) -- (1,0,1);
   \draw[] (0,0,1) -- (1,0,1) -- (1,1,1) -- (0,1,1) -- cycle;
\end{scope}
  
%normal cap
\begin{scope}
   [shift={(-.8,-.5)},x={(-.35cm,.35cm)},
    y={({cos(15)*6cm},{sin(15)*6cm})},
    z={({cos(91)*5cm},{sin(91)*5cm})},line join=round,fill opacity=0.5,black]
  \draw[] (0,0,0) -- (0,0,1) -- (0,1,1) -- (0,1,0) -- cycle;
   \draw[] (0,0,0) -- (1,0,0);
   \draw[] (1,0,0) -- (1,0,1);
   \draw[] (0,0,1) -- (1,0,1) -- (1,1,1) -- (0,1,1) -- cycle;
\end{scope}

\begin{scope}
   [shift={(-.8,0)},x={(-.35cm,.35cm)},
    y={({cos(6)*6cm},{sin(6)*6cm})},
    z={({cos(90)*5cm},{sin(90)*5cm})},line join=round,fill opacity=0.5,black]
  \draw[] (0,0,0) -- (0,0,1) -- (0,1,1) -- (0,1,0) -- cycle;
   \draw[] (0,0,0) -- (1,0,0);
   \draw[] (1,0,0) -- (1,0,1);
   \draw[] (0,0,1) -- (1,0,1) -- (1,1,1) -- (0,1,1) -- cycle;
\end{scope}

%medium plank
\begin{scope}
   [shift={(0,-.9)},x={(3cm,0cm)},
    y={({cos(60)*2.2cm},{sin(60)*2.2cm})},
    z={({cos(90)*5cm},{sin(90)*5cm})},line join=round,fill opacity=0.5, gray, line width=0.5mm]
  \draw[] (0,0,0) -- (0,0,1) -- (0,1,1);
   \draw[] (0,0,0) -- (1,0,0) -- (1,1,0);

  \draw[] (1,0,0) -- (1,0,1) -- (1,1,1) -- (1,1,0) -- cycle;
 \draw[] (0,0,1) -- (1,0,1) -- (1,1,1) -- (0,1,1) -- cycle;

\end{scope}

\draw[black, ->] (10,0) -- (11,0) node[right]{$\boldsymbol e_2(\theta)$};
\draw[black, ->] (10,0) -- (9.5,-0.8) node[below]{$\boldsymbol e_1(\theta)$};
\draw[black, ->] (10,0) -- (10,1) node[above]{$\boldsymbol e_3(\theta)$};

\node[] at (6,1)
{$P_{\theta}$};
\node[] at (2,2)
{$S_\tau$};
\node[] at (4,3)
{$R_\si$};

%\node[scale=0.8] at (1.38,4.7) {$\theta_{\tau_2,high}^+$};
% \node[scale=0.8] at (.8,5.2) {$\theta_{\tau_1,high}^+$};
% \node at (.8,6) {$\omega$}; 
 
\end{tikzpicture}
\caption{Relation between planks}
\label{relationbetweenplank}
\end{figure}

It remains to prove Lemma \ref{disjointlemma}. Before proving the lemma, we give some intuition on why the lemma should be true. See Figure \eqref{relationbetweenplank}. We first cover $N_{\lambda^2}\Ga_{K^{-1}}$ by gray planks $R_\si$ of dimensions $\sim\lambda^2\times\lambda\times 1$. Fix a $R_\si$, we draw all the black slabs $P_\theta$ of dimensions $\de\times 1\times 1$ whose corresponding $Q_\theta$ is contained in $R_\si$. Morally speaking, $P_\theta\cap R_\si\approx Q_\theta$ which is a $\de\times \lambda\times 1$-plank. One tricky thing is that different $P_\theta$ may have essentially same $Q_\theta$, which is the reason to introduce $S_\tau$ (the thick-black planks in the Figure of dimensions $\de\times\lambda\times 1$). Suppose we have a partition $R_\si=\bigcup_{\tau\prec\si}S_\tau$. We see that each $P_\theta\cap R_\si$ is contained in one of the $S_\tau$. If so, then we define $\theta\prec\tau.$ 
We can talk about the intuition on the numerology of these planks. 
\begin{enumerate}
    \item $\#\{R_\si\}=\lambda^{-1}$, 
    \item $\#\{Q_\theta:Q_\theta\subset R_\si\}\sim  \#\{\theta\prec\si\}=\de^{-1}\lambda$,\label{22}
    \item $\#\{S_\tau:S_\tau\subset R_\si\}\sim \frac{|R_\si|}{|S_\tau|}=\de^{-1}\lambda^2$,
    \item $\#\{Q_\theta:Q_\theta\subset S_\tau\}\sim \frac{\#\{Q_\theta:Q_\theta\subset R_\si\}}{\#\{S_\tau:S_\tau\subset R_\si\}}=\lambda^{-1}$. \label{44}
\end{enumerate}
By \eqref{22}, we see Property \eqref{property1} in Lemma \ref{disjointlemma} should be true. By \eqref{44}, we see Property \eqref{property2} and \eqref{property3} in Lemma \ref{disjointlemma} should be true.

 \begin{sloppypar}  \begin{proof}[Proof of Lemma \ref{disjointlemma}]
Recall $\be_1(\theta)=\ga(\theta),\be_2(\theta)=\ga'(\theta), \be_3(\theta)=\ga(\theta)\times\ga'(\theta)$. Defining 
$$\kappa(\theta)=\langle \be_2'(\theta),\be_3(\theta) \rangle (\gtrsim 1),$$
we have
\begin{align}
&\be_1'(\theta)=\be_2(\theta),\\
&\be_2'(\theta)=-\be_1(\theta)+\kappa(\theta)\be_3(\theta),\\
&\be_3'(\theta)=-\kappa(\theta)\be_2(\theta).
\end{align}

To prove Property \eqref{property1}, write $\theta=\si+\Delta$ with $|\Delta|\lesssim\lambda$. Any $\xi\in Q_\theta$ can be written as $\xi=a\be_1(\theta)+b\be_2(\theta)+c\be_3(\theta)$ with $|a|\le \de, |b|\le\lambda,|c|\le 1$. By Taylor's expansion, we have
\begin{align}
    \label{alsoprove}\xi=&a(\be_1(\si)+\Delta\be_2(\si))+ b(\be_2(\si)-\Delta\be_1(\si)+\Delta\kappa(\si)\be_3(\si))\\
    \nonumber&+c(\be_3(\si)-\Delta\kappa(\si)\be_2(\si))+O(\Delta^2)\\
    \nonumber=&(a-b\Delta)\be_1(\si)+(a\Delta+b-c\Delta\kappa(\si))\be_2(\si)+(b\Delta\kappa(\si)+c)\be_3(\si)+O(\Delta^2).
\end{align} 
One can easily check $\xi\in R_\si$.
The Property 
\eqref{property2} can also be proved by using \eqref{alsoprove}.
For the Property \eqref{property3}, we have proved a special case $\lambda\sim 1$ in the High case, but here we
need to do more work. We may assume $\lambda<<1$. Consider the plane
$$ \Pi_\theta=\{\xi_2\be_2(\theta)+\xi_3\be_3(\theta):\frac{1}{2}\lambda\le\xi_2\le\lambda, K^{-1}\le\xi_3\le 1\}. $$
We see $Q_\theta$ is the $\de$-neighborhood of $\Pi_\theta$. We just need to show: For any $\theta\in[0,1]$ and $\theta'=\theta+\Delta\in[0,1]$ with $C\lambda^{-1}\de\le \Delta \le C^{-1}$ (for some bounded $C$ to be determined later), if $\xi\in\Pi_\theta$, then
$$ \dist(\xi,\Pi_{\theta+\Delta})>10\de. $$
Write $\xi=a\be_2(\theta)+b\be_3(\theta)$, where $a\in[\frac{1}{2}\lambda,\lambda]$ and $b\in[K^{-1},1]$. 

We consider two scenarios: (1) $C\lambda^{-1}\de\le\Delta\le C\lambda$, (2) $C\lambda\le \Delta\le C^{-1}$.
If we are in the first scenario,
since the normal direction of $\Pi_{\theta+\Delta}$ is $\be_1(\theta+\Delta)$, it suffices to prove
\begin{equation}\label{disjoint2}
    |\be_1(\theta+\Delta)\cdot\big( a\be_2(\theta)+b\be_3(\theta) \big)|\ge 10\de. 
\end{equation} 
By Taylor's expansion, we have $\be_1(\theta+\Delta)=\be_1(\theta)+\Delta\be_2(\theta)+\frac{\Delta^2}{2}(-\be_1(\theta)+\kappa(\theta)\be_3(\theta))+o(\Delta^2)$. We see the left hand side of \eqref{disjoint2} is $| a\Delta+\frac{\Delta^2}{2}\kappa(\theta)b+o(\Delta^2) |\ge a\Delta-o(\Delta^2)=(a-o(\Delta))\Delta  \ge \frac{10}{C}a\Delta\gtrsim 10\de$, when $\Delta\le C\lambda$ and $\lambda<<1$. If we are in the second scenario, we show that 
\begin{equation}\label{disjoint3}
    \be_2(\theta+\Delta)\cdot (a\be_2(\theta)+b\be_3(\theta))\ge 10\lambda.
\end{equation}
By Taylor's expansion, we have
$\be_2(\theta+\Delta)=\be_2(\theta)+\Delta(-\be_1(\theta)+\kappa(\theta)\be_3(\theta))+O(\Delta^2)$.
We see the left hand side of \eqref{disjoint3}
is $|a+\Delta\kappa(\theta)b+O(\Delta^2)|\ge \Delta\kappa(\theta)b-a-O(\Delta^2)\ge 10\lambda$ if the constant $C$ is big enough.
\end{proof} \end{sloppypar}

% \section{Introduction} \label{areaproof}

\section{Proof of Theorem~\ref{thmpositive}}

For a small positive number $\delta$ and $E\subset [0, 1]$, we use $\Lambda_{\delta}(E)$ to denote a maximal $\de$-separated subset of $E$. By definition, $\#\Lambda_\de(E)\sim |E|_\de$. If $E=[0,1]$, then we abbreviate $\Lambda_{\delta}(E)$ as $\Lambda_\de$ and just choose it to be the $\de$-lattice points in $[0,1]$. 
A rectangular box of dimensions $\delta\times \delta\times 1$ will be referred to as a $\delta$-tube. For each $\theta\in \Lambda_{\delta}$, there is a set of finitely overlapping collection of $\delta$-tubes that cover $\R^3$ whose long sides are parallel to $\gamma(\theta)$. 

% We will use $\W_\theta$ to denote some subset of these tubes, where the subscript $\theta$ indicates the direction of the tubes. 

In order to prove Theorem \ref{thmpositive}, we need the following result about incidence estimate. 

\subsection{An incidence estimate}

\begin{theorem}\label{thmincidence}
Let $\Lambda_{\delta}$ be a $\delta$-net of $[0, 1]$ for some $\delta>0$. Given a small constant $\varepsilon>0$ and $0<\alpha\le 2$, let $\mu$ be a finite nonzero Borel measure supported in the unit ball in $\R^3$ with $c_\alpha(\mu):=\sup_{x\in\R^3,r>0}\frac{\mu(B(x,r))}{r^\alpha}<\infty$.  Suppose that $\W$ is a set of $\de$-tubes, with directions in $\{\ga(\theta):\theta\in\Lambda_{\delta}\}$, such that each $\W_{\theta}$ is disjoint, where we use $\W_{\theta}$ to denote  
the subset of tubes in $\W$ that points to direction $\ga(\theta)$. Suppose also that
\begin{equation}\label{220527e2_1a}
    \sum_{T \in \mathbb{W}} \chi_{T}(x) \gtrsim \delta^{\varepsilon-1}, \qquad \forall \, x \in \supp (\mu). 
\end{equation}
Then 
\begin{equation}
    |\W|\ge C_{\varepsilon,  \alpha} \cdot \mu(\R^3) c_{\alpha}(\mu)^{-1} \delta^{-(1+\alpha-O(\sqrt{\varepsilon}))},
\end{equation}
where the constant $C_{\varepsilon,  \alpha}$ is  allowed to depend on $ \alpha$ and $\varepsilon$ but not on $\delta$, and $O(\sqrt{\varepsilon})$ can be taken to be $10^{10}\sqrt{\varepsilon}$. We also remark that %we can assume $C_{\e,\alpha}\le C_{\e',\alpha}$ if we replace $C_{\e,\alpha}$ by $\min_{\e'\le \e}C_{\e,\alpha}$.
$C_{\varepsilon',\alpha}\le C_{\varepsilon'',\alpha}$ for $\varepsilon' \le \varepsilon''$.
\end{theorem}

\begin{proof}[Proof of Theorem \ref{thmincidence}] 
The main argument of the proof is similar to that of Theorem \ref{planediscrt}, except for the ``low case". 
%The proof is via an induction on $\varepsilon$. The base case of the induction is $10^{-2}\le \varepsilon\le 1$, which is trivial. To prove Theorem \ref{thmincidence} for $\varepsilon$, let us assume that we have proven it for $\widetilde{\varepsilon}\ge\frac{\varepsilon}{1-\sqrt{\varepsilon}}$. 

For a given $\delta$-tube $T_\theta\in \W_\theta$, we denote the dual slab of $T_\theta$ by $P_\theta$ which is of dimensions $\de^{-1}\times\de^{-1}\times 1$ and centered at the origin.
Recalling \eqref{coord}, we can write
$$ P_\theta=\left\{\sum_{i=1}^3 \xi_i \mathbf{e}_i(\theta) :|\xi_1|\le 1, |\xi_2|\le \de^{-1},|\xi_3|\le \de^{-1} \right\}. $$

\begin{remark}
{\rm The slab $P_\theta$ here is just the $\de^{-1}$-dilation of that in the proof of Theorem \ref{planediscrt}. So are the tubes $T_\theta$ and $P_{\theta,high}, P_{\theta,low}, P_{\lambda, \theta}$ that we will define right now.}
\end{remark}

Let $\phi_T: \R^3\to \R$ be a non-negative function with $\widehat{\phi_T}$ supported on $P_\theta$ and $\phi_T(x)\gtrsim 1$ for every $x\in T$. 
Set 
\begin{equation}
    f_\theta(x):=\sum_{T\in \W_\theta} \phi_T(x),\ \ \   f(x):=\sum_{T\in \W} \phi_T(x).
\end{equation} 
The assumption \eqref{220527e2_1a} implies that 
\begin{equation}
   f(x)\gtrsim \delta^{\varepsilon-1},
\end{equation}
for every $x\in \supp(\mu)$. 

Next we will do the frequency decomposition for $f_\theta$. Similar to Definition \ref{defdecomposition}, we make the following definitions.
\begin{definition} (See Figure \ref{decomposition}.)
Let $K=\de^{-\sqrt{\e}}$. Define the high part of $P_\theta$ as
$$ P_{\theta,high}:=\left\{\sum_{i=1}^3\xi_i\be_i(\theta): |\xi_1|\le 1, K^{-1}\de^{-1}\le|\xi_2|\le \de^{-1}, |\xi_3|\le \de^{-1} \right\}. $$
Define the low part of $P_\theta$ as 
$$ P_{\theta,low}:=\left\{\sum_{i=1}^3\xi_i\be_i(\theta): |\xi_1|\le 1, |\xi_2|\le K^{-1}\de^{-1}, |\xi_3|\le K^{-1}\de^{-1}\right\}. $$
For dyadic numbers $\lambda\in(\de^{1/2},K^{-1}]$, define 
$$ P_{\theta,\lambda}:=\left\{\sum_{i=1}^3\xi_i\be_i(\theta): |\xi_1|\le 1, \frac{1}{2}\lambda\de^{-1}\le|\xi_2|\le \lambda\de^{-1}, K^{-1}\de^{-1}\le|\xi_3|\le \de^{-1}\right\}. $$
In particular, we define
$$ P_{\theta,\de^{1/2}}:=\left\{\sum_{i=1}^3\xi_i\be_i(\theta): |\xi_1|\le 1, |\xi_2|\le \de^{-1/2},K^{-1}\de^{-1}\le|\xi_3|\le \de^{-1} \right\}. $$
\end{definition}
Similarly to \eqref{decompftheta} and \eqref{decompf}, we have
\begin{equation}
    f_\theta=f_{\theta,high}+f_{\theta,low}+ \sum_{\de^{1/2}\le\lambda\le K^{-1}}f_{\theta,\lambda},
\end{equation}
and
\begin{equation}
    f=f_{high}+f_{low}+ \sum_{\de^{1/2}\le\lambda\le K^{-1}}f_{\lambda}.
\end{equation}

Since $f(x)\gtrsim \de^{\e-1}$, there are two cases:\\

\noindent\fbox{Case 1: High case} We can find a Borel set $F$ satisfying 
\[ \mu(F) \gtrsim  \mu(\mathbb{R}^3), \]
and 
\begin{equation}\label{local2} \delta^{\varepsilon-1} \lesssim  \left\lvert f_{high}(x)+\sum_\lambda f_\lambda(x) \right\rvert, \quad \forall \, x \in F.  \end{equation}

\noindent\fbox{Case 2: Low case} We can find a Borel set $F$ satisfying
\[ \mu(F) \gtrsim \mu(\mathbb{R}^3), \]
and 
\begin{equation} \label{local} \delta^{\varepsilon-1} \lesssim  \left\lvert f_{low}(x) \right\rvert, \quad \forall \, x \in F. \end{equation}

Assume first that we are in the high case. We raise both sides of \eqref{local2} to the sixth power, integrate with respect to $d\mu$, and obtain 
\begin{equation}
    \mu(\R^3) \delta^{-6(1-\varepsilon)} \lesssim (\log\de^{-1})^{O(1)}\int |f_{high}|^6+\sum_\lambda |f_\lambda|^6 d\mu.
\end{equation}
Since the functions on the right hand side are locally constant on $\de$-balls, together with the upper density condition on $\mu$, we obtain 
\begin{equation}\label{220527e2_25} 
\mu(\mathbb{R}^3) c_{\alpha}(\mu)^{-1} \delta^{-6(1-\varepsilon) + 3-\alpha} \lesssim (\log\de^{-1})^{O(1)}\int |f_{high}(x)|^6+\sum_\lambda |f_\lambda(x)|^6 \, dx.
\end{equation}
We can just use \eqref{highpart} and \eqref{lambdapart} with $t=1$ and $\#\Theta\sim \de^{-1}$, noting there is a scaling difference, so that the right hand side above is bounded by
$$ \de^{-O(\sqrt{\e})}\de^{-2}|\W|. $$
It follows that 
\[ |\mathbb{W}| \gtrsim \mu(\mathbb{R}^3) c_\alpha(\mu)^{-1}  \delta^{-(1+\alpha-O(\sqrt{\varepsilon}))}.\]
This finishes the proof if we are in the high case \eqref{local2}.\\

Now we assume that we are in the low case \eqref{local}. For each $T\in \W$, the support of $\phi_T* \widecheck{\eta_{low}}$ is essentially a thickened tube of $T$ with dimensions $1\times K\delta\times K \delta$ where $K=\de^{-\sqrt{\e}}$. We use $\widetilde{T}$ to denote this thickened tube. Let $\widetilde{\W}$ be the collection of these thickened tubes obtained from $\W$, and we only keep those essentially distinct tubes (each tube intersects $\lesssim 1$ other tubes of those whose angle is within $\lesssim K \delta$ of its own, and every $T \in \mathbb{W}$ is contained in some $\widetilde{T}$ from $\widetilde{\W}$). 

% To be more precise, for two tubes $T_1, T_2$, let $\widetilde{T}_1, \widetilde{T}_2$ denote the enlarged tubes as above. We say that 
% \begin{equation}
%     \widetilde{T}_1\sim \widetilde{T}_2
% \end{equation}
% if 
% \begin{equation}
%     \widetilde{T}_1\cap \widetilde{T}_2\neq \emptyset, 
% \end{equation}
% and the angle between $T_1, T_2$ is $\lesssim K\delta$. Then $\widetilde{\W}$ is a collection of enlarged tubes $\widetilde{T}$ satisfying that
% \begin{equation}
%      \widetilde{T}_1\not\sim \widetilde{T}_2, \ \ \forall \widetilde{T}_1, \widetilde{T}_2\in \widetilde{\W}, \widetilde{T}_1\neq \widetilde{T}_2,
% \end{equation}
% and for every $T_1\in \W$, there exists $T_2$ with $\widetilde{T}_2\in \widetilde{\W}$ such that $\widetilde{T}_1\sim_{\rho} \widetilde{T}_2$.
Write $\widetilde{\W}$ as a disjoint union 
\begin{equation}
    \widetilde{\W}=\widetilde{\W}_{\heavy}\bigcup \widetilde{\W}_{\light}
\end{equation}
where $\widetilde{\W}_{\light}$ is the collection of thickened tubes $\widetilde{T}$ that contain $\le C^{-1}K^{3-\sqrt{\varepsilon}}$ tubes from $\W$. Here $C$ is a large universal constant which is much larger than the implicit (universal) constant in \eqref{local}. Note that 
\begin{align}
    \Norm{\sum_{\widetilde{T}\in \widetilde{\W}_{\light}} \sum_{\substack{T'\in \W: T'\subset \widetilde{T}}} \phi_{T'}* \widecheck{\eta_{low}}
    }_{\infty}
    &\le 
    \sum_{\widetilde{T}\in \widetilde{\W}_{\light}}
    \sum_{\substack{T'\in \W: T'\subset \widetilde{T}}} \Norm{  \phi_{T'}* \widecheck{\eta_{low}}
    }_{\infty}\\
    &\lesssim C^{-1} (\delta^{-1}K^{-1}) K^{-2} K^{3-\sqrt{\varepsilon}}=C^{-1}\de^{-1+\e}  
\end{align}
As a consequence, we see that \eqref{local} can be upgraded to \begin{equation}\label{220530e2_33}
    \delta^{\varepsilon-1}
    \lesssim \anorm{\sum_{\widetilde{T}\in \widetilde{\W}_{\heavy}} \sum_{\substack{T'\in \W: T'\subset \widetilde{T}}} \phi_{T'}* \widecheck{\eta_{low}}(x)}, \ \forall x\in F. 
\end{equation}
Next, by \eqref{220530e2_33} and the fact that if $x\in\widetilde{T}$ then $\#\{T\in\W:x\in T\subset \wt T\}\lesssim K$, we conclude that 
\begin{equation}
    \sum_{\widetilde{T}\in \widetilde{\W}_{\heavy}} \chi_{\widetilde{T}}(x)\gtrsim K^{-1}\delta^{\varepsilon-1}, \ \forall x\in F. 
\end{equation}
Write 
\begin{equation}
    \delta^{-1+\varepsilon}K^{-1}=(\delta K)^{-1+\widetilde{\varepsilon}}, \ \widetilde{\varepsilon}:=\frac{\varepsilon}{1-\sqrt{\varepsilon}}. 
\end{equation}
The tubes in $\widetilde{\W}_{\heavy}$ satisfy the induction hypothesis at the scale $K\de$ with the new parameter $\widetilde{\varepsilon}$. Hence 
\begin{align*} 
|\mathbb{W} | 
\geq 
|\widetilde{\mathbb{W}}_{\heavy} | K^{3-\sqrt{\varepsilon}} 
& 
\geq
%|\log \delta|^{-1}
\mu(\mathbb{R}^3) c_{\alpha}(\mu)^{-1} C_{\widetilde{\varepsilon},  \alpha} \left(\delta^{-1}K^{-1}\right)^{(\alpha+1- 10^{10}\sqrt{\widetilde{\varepsilon}} ) } K^{3-\sqrt{\varepsilon}}
%\\
% &
% \geq 
% \mu(\mathbb{R}^3) 
% c_{t}(\mu)^{-1} 
% C_{\varepsilon} 
% C_{\gamma} \left(
% \delta^{-1}\rho\right)^{(t+1- 10^{10}\sqrt{\widetilde{\epsilon}} ) } \rho^{-3+\sqrt{\varepsilon}} 
% &&\text{(monotonicity of $C_{\epsilon}$)}.
\end{align*}
Elementary computation shows that 
\begin{equation}
    \sqrt{\widetilde{\varepsilon}}-\sqrt{\varepsilon}=\frac{\widetilde{\varepsilon}-\varepsilon}{\sqrt{\widetilde{\varepsilon}}+\sqrt{\varepsilon}}\le \frac{\varepsilon}{2(1-\sqrt{\varepsilon})}\le \frac{3\varepsilon}{4}. 
\end{equation}
Therefore, 
\begin{equation}
\begin{split}
        |\W|& \ge C_{\widetilde{\varepsilon},  \alpha} \cdot \mu(\R^3) c_\alpha(\mu)^{-1}  \delta^{-(\alpha+1-10^{10}\sqrt{\varepsilon})} %|\log \delta|^{-1} 
        \delta^{10^{10}\frac{3\varepsilon}{4}} K^{3-(1+\alpha-10^{10}\sqrt{\widetilde{\varepsilon}})-\sqrt{\varepsilon}}\\
        & \ge C_{\widetilde{\varepsilon},  \alpha} \cdot \mu(\R^3) c_\alpha(\mu)^{-1}  \delta^{-(\alpha+1-10^{10}\sqrt{\varepsilon})} K^{2-\alpha} %|\log \delta|^{-1} 
        \delta^{10^{10}\frac{3\varepsilon}{4}} \delta^{-10^{10}\varepsilon+\varepsilon}.
\end{split}
\end{equation}
Since $\alpha \le 2$, the induction closes. \end{proof}

We have the following corollary of Theorem \ref{thmincidence}.

\begin{corollary}\label{220514corollary2_2}
Let $\alpha\in (0, 3)$ and $\alpha_1\in (0, \min\{2, \alpha\})$. Let $\mu$ be a finite non-zero Borel measure supported on the unit ball in $\R^3$ with $c_{\alpha}(\mu)\le 1$. Let $\delta>0$ be a small number. Let $\Lambda_{\delta}$ be a $\delta$-net of $[0,1]$. For each $\theta\in \Lambda_{\delta}$, let $\D_{\theta}$ be a disjoint collection of at most $\mu(\R^3)\delta^{-\alpha_1}$ balls of radius $\delta$ in $\pi_{\theta}(\R^3)$. Then there exists $\varepsilon>0$, depending only on $\alpha$ and $\alpha_1$, such that 
\begin{equation}\label{220515e2_3}
    \delta \sum_{\theta\in \Lambda_{\delta}}  (\pisharp\mu)\pnorm{\bigcup_{D\in \D_{\theta}}D} \le C(\alpha, \alpha_1)\mu(\R^3) \delta^{\varepsilon}, 
\end{equation}
where the constant $C(\alpha, \alpha_1)$ depends on $\alpha, \alpha_1$, but not on $\delta$. 
\end{corollary}
\begin{proof}[Proof of Corollary \ref{220514corollary2_2}]
We argue by contradiction. Suppose that for every $\varepsilon>0$, there exist $\delta>0$ and $\D_{\theta}$, a disjoint collection of at most $\mu(\R^3)\delta^{-\alpha_1}$ balls of radius $\delta$ in $\pi_{\theta}(\R^3)$ for each $\theta\in \Lambda_{\delta}$, such that \eqref{220515e2_3} fails. Note that for each $D_\theta\in\D_\theta$, $\pi_\theta^{-1}(D_\theta)\cap B^3(0,1)$ is a $\de$-tube. We denote these tubes by 
\begin{equation}
    \W_{\theta}:=\{\pi_\theta^{-1}(D_\theta)\cap B^3(0,1): D_{\theta}\in \D_{\theta}\},
\end{equation}
and write $\W:=\cup_{\theta} \W_{\theta}$. Define 
\begin{equation}
    F:=\Big\{x\in \supp(\mu): \delta \sum_{\theta\in \Lambda_{\delta}} \sum_{T\in \W_{\theta}} \chi_T(x)\ge C( \alpha, \alpha_1) \delta^{\varepsilon}/2 \Big\}
\end{equation}
Note that by our contradiction assumption, we have 
\begin{equation}
   \begin{split}
         C( \alpha, \alpha_1) \mu(\R^3) \delta^{\varepsilon}& \le \delta \sum_{\theta\in \Lambda_{\delta}}  (\pisharp\mu)\pnorm{\bigcup_{D\in \D_{\theta}}D}\\
         & = \delta  \int_{\R^3}  \sum_{\theta}\sum_{T\in \W_{\theta}} \chi_T(x) d\mu(x).
   \end{split}
\end{equation}
This further implies 
\begin{equation}
    \mu(F)\gtrsim C( \alpha, \alpha_1)\mu(\R^3)\delta^{\varepsilon}. 
\end{equation}
Note that by definition, for every $x\in F$, it holds that \begin{equation}
    \sum_{T\in \W} \chi_T(x) \gtrsim \delta^{-1+\varepsilon}. 
\end{equation}
We apply Theorem \ref{thmincidence} to the measure $\mu$ restricted to $F$ and obtain 
\begin{equation}
    |\W|\gtrsim \delta^{\varepsilon} \delta^{-1-\min\{2, \alpha\}+O(\sqrt{\varepsilon})}\mu(\R^3).
\end{equation}
By pigeonholing, this further implies that there exists $\theta$ such that 
\begin{equation}
    |\W_{\theta}|\gtrsim \delta^{\varepsilon} \delta^{-\min\{2, \alpha\}+O(\sqrt{\varepsilon})}\mu(\R^3).
\end{equation}
This is a contradiction to the assumption that $|\D_{\theta}|\lesssim \delta^{-\alpha_1}\mu(\R^3)$ when $\varepsilon$ is chosen to be small enough. 
% Let $\Lambda'_{\delta}$ denote the collection of all $\theta\in \Lambda_{\delta}$ such that 
% \begin{equation}
%     (\pisharp\mu)\pnorm{\bigcup_{D\in \D_{\theta}}D}\ge C( \alpha, \alpha^*) \delta^{\varepsilon}/2.
% \end{equation}
% Note that by assumption each term under the sum in \eqref{220515e2_3} is $\le 1$.
% Consequently $|\Lambda'_{\delta}|\gtrsim \delta^{-1+\varepsilon}$.
\end{proof}

The corollary below is a special case of Corollary~\ref{220514corollary2_2}. It is recorded below for a later application.
\begin{corollary}
\label{alphastars} Let $\alpha \in [2,3]$ and $\alpha^* \in \left[0, 2\right)$. Let $\mu$ be a finite non-zero Borel measure supported on the unit ball in $\R^3$ with $c_{\alpha}(\mu)\le 1$. Let $\delta>0$ be a small number. Let $\Lambda_{\sqrt{\delta}}$ be a $\sqrt{\delta}$-net of $[0, 1]$.  For each $\theta\in \Lambda_{\sqrt{\delta}}$, let $\D_{\theta}$ be a disjoint collection of at most $\delta^{-\frac{\alpha^*+1}{2}}\mu(\R^3)$ rectangles of dimension $\delta\times \sqrt{\delta}$ in $\pi_{\theta}(\R^3)$ whose long sides point in the $\gamma'(\theta)$ direction. Then there exists $\varepsilon>0$, depending only on $\alpha$ and $\alpha^*$, such that 
\begin{equation}
    \sqrt{\delta} \sum_{\theta\in \Lambda_{\sqrt{\delta}}}  (\pisharp\mu)\pnorm{\bigcup_{D\in \D_{\theta}}D} \le C( \alpha, \alpha^*) \delta^{\varepsilon}\mu(\R^3), 
\end{equation}
where the constant $C( \alpha, \alpha^*)$ depends on $\gamma, \alpha, \alpha^*$, but not on $\delta$. 
\end{corollary}
\begin{proof}[Proof of Corollary \ref{alphastars}]
For each $\theta\in \Lambda_{\sqrt{\de}}$ and recalling \eqref{coord}, define
\begin{equation}
    U_{\theta}:=\left\{ \sum_{i=1}^3 x_i \be_i(\theta): |x_1|\le 1, |x_2|\le \sqrt\de, |x_3|\le \de \right\},
\end{equation}
which is a $\de\times\sqrt\de\times 1$-plank.
By a simple geometric observation, we have that 
\begin{equation}
    \pi_{\theta'}(W_{\theta})\subset C \pi_{\theta}(W_{\theta}),
\end{equation}
for every $\theta,\theta'\in [0, 1]$ with $|\theta'-\theta|\le \sqrt{\delta}$, where $C$ is a constant depending only on $\gamma.$ 
For each $\theta\in\Lambda_{\sqrt\de}$ we have a set of $\de\times\sqrt\de\times 1$-planks 
$$\U_\theta=\{ \pi_\theta^{-1}(D_\theta)\cap B^3(0,1) \}_{D_\theta\in\D_\theta},$$ 
which are essentially the translates of $U_\theta$. For each $\theta'\in\Lambda_\de$, we choose a $\theta\in\Lambda_{\sqrt\de}$ with $|\de-\de'|\le \sqrt{\de}$. We partition each plank in $\U_\theta$ into $\de\times \de\times 1$-tubes with direction $\ga(\theta')$ and denote all these tubes by $\W_{\theta'}$. We also define $\D_{\theta'}'=\pi_\theta(\W_{\theta'})$ which is a collection of at most $\mu(\R^3)\de^{-\frac{\alpha*}{2}-1}$ balls of radius $\de$ in $\pi_\theta(\R^3)$.

Now we can apply Corollary \ref{220514corollary2_2} to the sets $\{\D_{\theta'}'\}_{\theta'\in\Lambda_\de}$ with $\alpha_1=\frac{\alpha^*}{2}+1$.
We obtain
$$ \sqrt{\delta} \sum_{\theta\in \Lambda_{\sqrt{\delta}}}  (\pisharp\mu)\pnorm{\bigcup_{D\in \D_{\theta}}D}\sim \delta \sum_{\theta\in \Lambda_{\delta}}  (\pisharp\mu)\pnorm{\bigcup_{D\in \D_{\theta}'}D} \le C(\alpha, \alpha^*)\mu(\R^3) \delta^{\varepsilon}. $$ \end{proof}

%\medskip

% \begin{corollary} Let $\gamma: [a,b] \to\ZS^2$ be a $C^2$ curve with $\det(\gamma, \gamma', \gamma'') \neq 0$. If $A \subseteq \mathbb{R}^3$ is a Borel set with $\dim A \leq 2$, then $\dim \pi_{\theta}(A) = \dim A$ for a.e.~$\theta \in [a,b]$. \end{corollary}

%\medskip

By Frostman's Lemma (see for instance \cite[page 112]{mattila}), Theorem \ref{thmpositive} is an immediate consequence of the following theorem. 

\begin{theorem}\label{thmabsconti}
Let $\gamma: [0, 1] \to\ZS^2$ be $C^2$ and non-degenerate. If $\mu$ is a compactly supported Borel measure on $\mathbb{R}^3$ such that $c_{\alpha}(\mu) < \infty$ for some $\alpha >2$, then $\pi_{\theta\#} \mu$ is absolutely continuous with respect to $\mathcal{H}^2$, for a.e.~$ \theta \in [0, 1]$. 
\end{theorem}

To prove Theorem \ref{thmabsconti}, we will cut $\gamma$ into finitely many pieces with a number depending only on $\gamma$, and work on one piece. From now on, we assume that $\gamma: [0, \aaa]\to\ZS^2$ is $C^2$ and non-degenerate, and satisfy 
\begin{equation}
    \gamma(0)=(0, 0, 1), \ \ \gamma'(0)=(1, 0, 0), \ |\gamma'(\theta)|=1, \forall \theta\in [0, \aaa]. 
\end{equation}
Here $\aaa>0$ is sufficiently small depending on $\gamma$.\\

\subsection{Decomposition of the frequency space}
In this subsection, we discuss the decomposition of the frequency space $\R^3_\xi$. Recall the cone that we are considering:
$$\Ga=\{r(\ga\times\ga')(\theta):r\in\R, \theta\in[0,a]\}.$$
For any $\xi=r(\ga\times\ga')(\theta)\in\Ga$, there are three directions that we would like to specify: the normal direction $\ga(\theta)$; the tangent direction $\ga'(\theta)$; and the flat (or radial) direction $\ga\times\ga'(\theta)$. We want to decompose $\R^3_\xi$ into regions according to the distance from the origin, the distance from the cone $\Ga$, and the angular parameter. We give the precise definition below.

For a given integer $k$, define 
\begin{equation}
    \Theta_k:=\big(2^{-k/2}\N\big)\cap [0, a]. 
\end{equation}
   
Fix $k\le j$. If $k<j$ and $\Box \in \{+,-\}$, then we define the plank in the forward/backward light cone
\begin{multline}\label{taudefn} \tau(\theta,j,k, \Box) =  \Big\{ \lambda_1 \left(\gamma \times \gamma'\right)(\theta) + \lambda_2 \gamma'(\theta) + \lambda_3 \gamma(\theta) \\
:
 |\lambda_1| \sim 2^{j}, |\lambda_2| \le C_{\gamma}^{-1}2^{-k/2 + j}, |\lambda_3| \sim 2^{-k + j} , \sgn \lambda_1 = \Box = -\sgn \lambda_3\Big\}. 
\end{multline}
Here $C_{\gamma}>0$ is some large constant that depends only on $\gamma$. It is chosen such that the distance from $\tau(\theta, j, k)$ to cone $\{(\gamma\times \gamma')(\theta): \theta\in [0, a]\}$ is comparable to $2^{j-k}$. Let us digest a little bit about the plank $\tau(\theta,j,k)$: $\lambda_1\sim 2^j$ is the distance from the origin; $|\lambda_3|\sim 2^{-k+j}$ is the distance from the cone $\Ga$; and $|\lambda_2|\lesssim 2^{-k/2+j}$ is the angular parameter of the plank.
 
\begin{sloppypar} If $k=j$, then for $\Box \in \{+,-\}$ we define 
 \begin{multline} \label{taudefn3} \tau(\theta,j,j, \Box) = \\ \Big\{ |\lambda_1| \left(\gamma \times \gamma'\right)(\theta) + \lambda_2 \gamma'(\theta) + \lambda_3 \gamma (\theta) :
 \lambda_1 \sim 2^{j}, \quad |\lambda_2| \lesssim 2^{j/2}, \quad |\lambda_3| \lesssim 1, \sgn \lambda_1 = \Box \Big\}. 
 \end{multline}
 Let 
 \begin{equation}
     \Lambda_{j, k}:=\{\tau(\theta, j, k, \Box): \theta\in \Theta_k, \Box \in \{-,+\} \}, 
 \end{equation}
 and 
 \begin{equation}
     \Lambda_j:=\bigcup_{k\le j}\Lambda_{j, k}, \ \ \Lambda:=\bigcup_{j\in \N}\Lambda_j.  
 \end{equation}
Roughly speaking, for $k < j$, $\Lambda_{j,k}$ forms a canonical covering of the part of $\{\xi\in\R^3: |\xi|\sim 2^j; \dist(\xi,\Ga)\sim 2^{-k+j}\}$ outside the cone $\Gamma$ by $2^j\times 2^{-k/2+j}\times 2^{-k+j}$-planks. Each $\Lambda_j$ forms a covering of the $\sim 1$-neighbourhood of $\{\xi \in \Gamma:|\xi|\sim 2^j\}$. %$\Lambda$ forms a covering of $\{\xi:|\xi|\gtrsim 1\}$.  
\end{sloppypar}

For each $\tau\in \Lambda$, we define $\T_{\tau}$ to be a set of planks of dual dimensions to $\tau$ (but scaled by $2^{k\delta}$ in each direction where $\delta>0$ and $\tau \in \Lambda_{j,k}$) and forming a finitely overlapping covering of $\R^3$. We will refer to $\T_\tau$ as the wave packets determined by the plank $\tau$. 
Now, we discuss the wave packet decomposition. For each $\tau\in\Lambda$, we can choose a smooth bump function $\psi_\tau$ supported in $2\tau$ and choose a smooth bump function $\psi_0$ supported in the unit ball, so that we have the partition of the unity
$$\psi_0(\xi)+\sum_{\tau\in \Lambda}\psi_\tau(\xi)=1,%\ \ \ \textup{for~}\xi\in\R^3.
$$
on the union of the $\tau$'s. For each $T\in\T_\tau$, we can choose a smooth function $\eta_T$ which is essentially supported in $T$ (with rapidly decaying tail outside of $T$), such that $\supp\wh \eta_T\subset \tau$ and 
$$\sum_{T\in\T_\tau}\eta_T(x)=1.$$
%From now on, we also view $\mu$ as a function on $\R^3$ by smooth approximation. For any set $U\subset \R^3$, we use $\mu(U)$ to denote $\int_U |\mu|$.

For any $T\in\T_\tau$, we define the wave packet
$$ M_T\mu:=\eta_T\big( \mu*\widecheck{\psi_\tau} \big). $$

% Given a box $\tau \in \Lambda_{j,k}$ and $T \in \mathbb{T}_{\tau}$, define 
% \[ M_T \mu = \eta_T \left( \mu \ast \widecheck{\psi_{\tau} } \right), \ \  M_T \mu = \widecheck{\psi'_{\tau}}*\Big(\eta_T \left( \mu \ast \widecheck{\psi_{\tau} } \right)\Big)\]
% for each compactly supported finite Borel measure $\mu$.

%Let $\phi$ be a bump function equal to 1 on $B_3(0,1)$ which vanishes outside $B_3(0,2)$.
\begin{lemma}\label{pushtwo} 
For $\tau \in \Lambda_{j,k}$ and $T\in \T_{\tau}$, we have  
\[ \left\lVert M_T \mu \right\rVert_{L^1(\R^3)} \lesssim 2^{3k \delta}\mu(2T) + C_N 2^{-kN} \mu\left( \mathbb{R}^3\right), \]
for every $N \geq 1$. 
\end{lemma}
\begin{proof}[Proof of Lemma \ref{pushtwo}]
Note that $|\widecheck{\psi_\tau}(x)|\lesssim \phi_{\tau^*}(x)$, where $\phi_{\tau^*}(x)$ is an $L^1$ normalized function essentially supported in $\tau^*$ (the dual plank of $\tau$). So, we have $|\eta_T|*|\widecheck{\psi_\tau}|\lesssim |\eta_T|$. Therefore,
$$ \int |M_T\mu|\lesssim \int |\mu| \big(|\eta_T|*|\widecheck{\psi_\tau}|\big)\lesssim \int |\mu||\eta_T|\lesssim 2^{3k\de}\mu(2T)+C_N2^{-k N}\mu(\R^3).  $$
\end{proof}

\begin{lemma} \label{nonstat} 
For $\tau \in \Lambda_{j,k}$ and $\theta\in [0, a]$ with 
\begin{equation} \label{angleassumption2} | \theta - \theta_{\tau} | \geq 2^{-k(1/2 - \delta)},
\end{equation}
it holds that 
\begin{equation}
    \left\lVert \pi_{\theta \#} M_T f \right\rVert_{L^1(\mathcal{H}^2 )} \lesssim_{\delta, N} 2^{-kN} |\tau| \|f\|_{L^1(\R^3)}. 
\end{equation}
for every $N \geq 1$, $T \in \mathbb{T}_{\tau}$ and $f \in L^1(\mathbb{R}^3)$.
\end{lemma} 
\begin{proof}[Proof of Lemma \ref{nonstat}] 
We start by writing
\begin{equation}
    M_T f(x)=\eta_T(x)\int_{\R^3} \widehat{f}(\xi)\psi_{\tau}(\xi) e^{i\inn{x}{\xi}}d\xi.
\end{equation}
By identifying the complex measure $\pi_{\theta \#} M_Tf$ with its Radon-Nikodym derivative with respect to $\mathcal{H}^2$, we obtain 
\begin{equation}
    \begin{split}
        \pisharp M_T f(x)=\int_{\R^3} \widehat{f}(\xi)\psi_{\tau}(\xi) \Big[\int_{\R} \eta_T(x+t\gamma(\theta))  e^{it\inn{\gamma(\theta)}{\xi}} dt\Big] e^{i\inn{x}{\xi}}d\xi
    \end{split}
\end{equation}
for every $x \in \gamma(\theta)^{\perp}$. It therefore suffices to show that
\[ \left\lvert \int_{\mathbb{R} } \eta_T(x + t\gamma(\theta) ) e^{i t \langle \xi, \gamma(\theta) \rangle } \, dt \right\rvert \lesssim_N 2^{-kN}, \quad \forall \,  \xi \in \tau, \quad \forall \, x \in \mathbb{R}^3.  \] 
Integration by parts will finish the proof. For more details, we refer to Lemma 2.6 in \cite{harris2019improved}.
\end{proof}

\subsection{Proof of Theorem \ref{thmabsconti}: good part and bad part}
The main idea is to divide the wave packets into two parts, called the good part and the bad part. We will prove an $L^1$ estimate for the bad part and an $L^2$ estimate for the good part.
%Recalling the wave packet decomposition in the last subsection, we can write
%$$ \mu=\mu*\widecheck{\psi_0}+\sum_{\tau\in\Lambda}\mu*\widecheck{\psi_\tau}=\mu*\widecheck{\psi_0}+\sum_{\tau\in\Lambda}\sum_{T\in\T_\tau}\eta_T(\mu*\widecheck{\psi_\tau})=\mu*\widecheck{\psi_0}+\sum_{\tau\in\Lambda}\sum_{T\in\T_\tau}M_T\mu. $$

Let $\alpha>2$ be as in Theorem \ref{thmabsconti}. Let $\epsilon>0$ be a small number (note $\epsilon$ is different from $\e$) and let $\alpha_0>0$ be determined later (we will later let $\alpha_0\nearrow 2$). For $j\ge k$ and $\tau\in \Lambda_{j, k}$, define

\[ \mathbb{T}_{\tau,b} := \left\{ T \in \mathbb{T}_{\tau} : \mu(4T) \geq C 2^{- k(\alpha_0+1)/2-\alpha(j-k) + 10^3 j \epsilon} \right\}, \qquad \mathbb{T}_{\tau,g} = \mathbb{T}_{\tau} \setminus \mathbb{T}_{\tau,b}. \]
Define 
\begin{equation} \label{Linfty} 
\mu_b = \sum_{j\in \N} \sum_{\substack{k \in [j \epsilon, j]}} \sum_{\tau \in \Lambda_{j,k}} \sum_{T \in \mathbb{T}_{\tau,b}} M_T\mu,\ \ \ \mu_g=\mu-\mu_b.
\end{equation}
We remark that the wave packets of $\mu_b$ are those that have heavy $\mu$-mass and not too far away from the cone $\Ga$. 
We have

\[ \mu  = \mu_g + \mu_b. \]
We remark that 
$\mu_g = \mu_{g, \alpha, \epsilon, \alpha_0}$ and $\mu_b = \mu_{b, \alpha, \epsilon, \alpha_0}$ depends on parameters $\alpha,\epsilon,\alpha_0$, but for simplicity we just omit them.

Theorem \ref{thmabsconti} follows from Lemma \ref{badpart} and Lemma \ref{goodpart} below. \\

\begin{lemma} \label{badpart}
Let $\alpha>2$ and $\epsilon\ll 1$. Fix $\alpha_0<2$. For all Borel measures $\mu$ supported on the unit ball in $\R^3$ with $c_{\alpha}(\mu) \leq 1$, it holds that 
\[\int \int \left\lvert \pi_{\theta \#} \mu_{b} \right\rvert \, d\mathcal{H}^2 d\theta\lesssim 1, \]
where $\mu_b$ is defined by \eqref{Linfty}, and the implicit constant depends on $\alpha, \alpha_0, \epsilon$ and $\mu$. 
\end{lemma}

\begin{lemma} \label{goodpart} Let $\alpha>2$. Then for $\alpha_0<2$ sufficiently close to $2$, and $\epsilon>0$ small enough depending on $\alpha$ and $\alpha_0$, and $\delta \ll \epsilon$,
\begin{equation} \label{goal2} \int \int \left\lvert \pi_{\theta\#} \mu_g \right\rvert^2 \, d\mathcal{H}^2 \, d\theta %= \int_a^b \int \left\lvert \pi_{\theta\#} \mu_{g,j_0,J, \alpha, \epsilon, \delta, \alpha_0} \right\rvert^2 \, d\mathcal{H}^2 \, d\theta \\
\lesssim 1, 
\end{equation}
where the implicit constant depends on $\alpha, \alpha_0, \epsilon$ and $\mu$. 
\end{lemma} 

\medskip

%\subsubsection{Proof of Lemma \ref{badpart}}
\begin{proof}[Proof of Lemma \ref{badpart}]
By definition, we first write 
\begin{equation} \label{badbound} 
\int \int \left\lvert \pi_{\theta \#} \mu_b \right\rvert d\mathcal{H}^2 d\theta \leq \int \sum_{j\in \N} \sum_{\substack{ k \in [j\epsilon,j]} } \sum_{\tau \in \Lambda_{j,k}} \sum_{T \in \mathbb{T}_{\tau,b} }  \int \left\lvert \pi_{\theta \#} M_T\mu  \right\rvert d\mathcal{H}^2 d\theta .
\end{equation}
By the triangle inequality, this is 
\begin{align}
&\label{firstterm} \le  \sum_{j} \sum_{\substack{ k \in [j\epsilon,j]} }    \int \int \sum_{\substack{\tau \in \Lambda_{j,k}: \\
\left\lvert \theta_{\tau} - \theta \right\rvert < 2^{k(-1/2+\delta) }}} \sum_{T \in \mathbb{T}_{\tau,b} } \left\lvert \pi_{\theta \#} M_T\mu \right\rvert d\mathcal{H}^2 d\theta\\
\label{secondterm} &\qquad + \sum_{j} \sum_{\substack{ k \in [j\epsilon,j]} }    \int\int \sum_{\substack{\tau \in \Lambda_{j,k}: \\
\left\lvert \theta_{\tau} - \theta \right\rvert \geq 2^{k(-1/2+\delta) }}} \sum_{T \in \mathbb{T}_{\tau,b} } \left\lvert \pi_{\theta \#} M_T\mu \right\rvert d\mathcal{H}^2 d\theta. \end{align}
By Lemma~\ref{nonstat}, the contribution from \eqref{secondterm} is
\[ \lesssim_{\delta, \epsilon,N} \sum_j \sum_{k\in [\epsilon j, j]} 2^{3j-\frac{3}{2}k-kN} \mu\left( \mathbb{R}^3 \right)\lesssim_{\delta, \epsilon} \sum_j \sum_{k\in [\epsilon j, j]} 2^{-j} \mu\left( \mathbb{R}^3 \right) \lesssim \mu(\R^3), \]
By choosing $N>100\e^{-1}$.

To estimate \eqref{firstterm}, we discretize the integration in $\theta$ and bound it by 
\begin{equation}
    \sum_{j} \sum_{\substack{ k \in [j\epsilon,j] } }  \sum_{\theta\in \Theta_k}  2^{-k/2}\int \sum_{\substack{\tau \in \Lambda_{j,k}: \\
\left\lvert \theta_{\tau} - \theta \right\rvert < 2^{k(-1/2+\delta) }}} \sum_{T \in \mathbb{T}_{\tau,b} } \left\lvert \pi_{\theta \#} M_T\mu \right\rvert d\mathcal{H}^2
\end{equation}
By Lemma~\ref{pushtwo} the contribution from \eqref{firstterm} is
\begin{align} \notag &\lesssim \mu\left( \mathbb{R}^3 \right) +  \sum_{j} \sum_{\substack{ k \in [j\epsilon,j]} } \sum_{\theta} 2^{-k/2}\sum_{\substack{\tau \in \Lambda_{j,k}: \\
\left\lvert \theta_{\tau} - \theta \right\rvert < 2^{k(-1/2+\delta) }}} \sum_{T \in \mathbb{T}_{\tau,b} }  2^{3j\delta} \mu(2T)\\
\label{intermediate} &\lesssim \mu\left( \mathbb{R}^3 \right)+ \sum_{j} \sum_{\substack{ k \in [j\epsilon,j]} } \sum_{\theta}   2^{-k/2}2^{100j \delta}\mu (B_{j,k}(\theta)) , \end{align}
where 
\[ B_{j,k}(\theta) = \bigcup_{\substack{ \tau \in \Lambda_{j,k}: \\ |\theta_{\tau} - \theta | < 2^{k(-1/2+\delta) }}} \bigcup_{T \in \mathbb{T}_{\tau,b} } 2T. \] 
For fixed $j$ and $k$, let $\{B_l \}_l$ be a finitely overlapping cover of the unit ball in $\R^3$ by balls of radius $2^{-(j-k)}$. For each $\theta$ and $l$ let 
\[ B_{j,k,l}(\theta) = \bigcup_{\substack{ \tau \in \Lambda_{j,k}: \\ |\theta_{\tau} - \theta | <2^{k(-1/2+\delta) }}} \bigcup_{\substack{T \in \mathbb{T}_{\tau,b}:  \\ 2T \cap B_l \neq \emptyset }} 2T.\]
Let $\mu_{j,k}$ be the pushforward of $\mu$ under $x \mapsto 2^{j-k-2k\delta} x$. Denote 
\begin{equation}
    B'_{j, k, l}(\theta):=2^{j-k-2k\delta} \cdot B_{j,k,l}(\theta), \ \ \widetilde{B}_{l} = \left\{ 2^{j-k-2k\delta}b_l + y : |y| \leq 1 , y\in \R^3\right\},
\end{equation}
with $b_l$ the centre of $B_l$, and define
% \[ B_{j,k,l}' = \left\{( \theta, x) \in [a,b] \times \mathbb{R}^3 : x \in 2^{j-k-2k\delta} \cdot B_{j,k,l}(\theta) \right\}, \]
\begin{equation}
    \widetilde{\mu}_{j,k,l} =2^{\alpha(j-k-2k\delta)} \cdot \mu_{j,k}\chi_{\widetilde{B}_l}.
\end{equation}
Then 
\begin{equation}\label{tildemu}
    \sum_{\theta} \mu(B_{j, k}(\theta)) \le \sum_{\theta} \sum_{l} \mu(B_{j, k, l}(\theta))\le \sum_l 2^{-\alpha(j-k-2k\delta)} \sum_{\theta} \widetilde{\mu}_{j, k, l}(B'_{j, k, l}(\theta)).
\end{equation}
Note that for each $\theta$, the set $B'_{j, k, l}(\theta)$ is contained in a union of planks of dimensions $1\times 2^{-k/2}\times 2^{-k}$; the number of planks is
\begin{equation}
    \lesssim 2^{k\frac{\alpha_0+1}{2}+C\delta k} \widetilde{\mu}_{j, k, l}(\R^3),
\end{equation}
for some large constant $C$, and each plank overlaps $\lesssim 2^{10k\delta}$ of the others. Moreover $c_{\alpha}\left(\widetilde{\mu}_{j,k,l}  \right) \leq 1$ and $\widetilde{\mu}_{j,k,l} $ is supported in a ball of radius 1. Therefore, by applying the triangle inequality and Corollary \ref{alphastars}, we can find $\delta'$, depending only on $\alpha$ and $\alpha_0$, such that 
\begin{multline*}
    2^{-k/2}\sum_{\theta}\widetilde{\mu}_{j, k, l}(B'_{j, k, l}(\theta))\lesssim 2^{-C'\delta'k+C'\delta k}  \widetilde{\mu}_{j, k, l}(\R^3) \\ \lesssim 2^{-C'\delta'k+C'\delta k} 2^{\alpha(j-k-2k\delta)}\mu_{j, k}(\widetilde{B}_l),
\end{multline*}
for some large constant $C'$, whose precise value is not important. Putting this into \eqref{tildemu} yields
 \begin{equation}
     \eqref{tildemu}\lesssim 2^{k/2} 2^{-C'\delta'k+C'\delta k} \mu(\R^3). 
 \end{equation}
Substituting this into \eqref{intermediate} and then \eqref{badbound} gives
\begin{equation}
    \eqref{badbound}\lesssim \sum_{j}\sum_{\substack{k\in [\epsilon j, j]}} 2^{100 \delta j} 2^{-C'\delta'k+C'\delta k}\mu(\R^3)
\end{equation}
Recall that $\delta'$ depends only on $\alpha$ and $\alpha_0$. We just need to pick $\delta$ to be sufficiently small, and will finish the proof. \end{proof} 

%\bigskip

%\subsubsection{Proof of Lemma \ref{goodpart}}
\begin{proof}[Proof of Lemma \ref{goodpart}] Take $\epsilon \ll \min\{\alpha-2,2-\alpha_0\}$. Given $x\in \pi_{\theta}(\R^3)$, note that 
\begin{equation}\label{pimu}
    \pisharp \mu_g(x)=\int \mu_g(x+t\gamma(\theta))dt. 
\end{equation}
Fix the coordinate $(\be_1,\be_2,\be_3)=(\ga'(\theta),\ga(\theta)\times\ga'(\theta),\ga(\theta))$. Any $x\in\pi_\theta(\R^3)$ can be written in this coordinate as $x=(x_1,x_2,0)$. We can also rewrite \eqref{pimu} as 
$$ \pi_{\theta\#}\mu_g(x_1,x_2)=\int \mu_g(x_1,x_2,t)dt. $$
Doing the Fourier transfom in the $(\be_1,\be_2)$-plane, we have \begin{multline*} (\pi_{\theta\#}\mu_g)^\wedge(\eta_1,\eta_2)=\int \mu_g(x_1,x_2,t)e^{-i(x_1\eta_1+x_2\eta_2)}  \, dt \\
=  \wh{\mu_g}(\eta_1,\eta_2,0)=\wh{\mu_g}(\eta_1 \gamma'(\theta)+ \eta_2 ( \gamma \times \gamma')(\theta)). \end{multline*}
By Plancherel's theorem, 
\begin{equation} \label{localise} \int \int \left\lvert \pi_{\theta \#} \mu_g \right\rvert^2 \, d\mathcal{H}^2 \, d\theta
= \int \int_{\mathbb{R}^2} \left\lvert \widehat{\mu_g}\left( \eta_1 \gamma'(\theta)+ \eta_2 \left( \gamma \times \gamma'\right)(\theta) \right)\right\rvert^2 \, d\eta \, d\theta. \end{equation}
Roughly speaking, 
\[  \mu_g = \mu_0+\mu_{g, 1}+\mu_{g, 2},\]
where $\mu_0$ is roughly $\mu(\R^3)\Id_{B^3(0,1)}$ with rapidly decaying tail outside $B^3(0,1)$, $\mu_{g,1}$ is the sum of good wave packets which have controlled mass, and $\mu_{g,2}$ is the sum of wave packets which are far away from the cone $\Ga$. A formula for $\mu_{g,1}$ is
%\begin{equation}
%    \begin{split}
%        \mu_g=& \mu*\widecheck{\psi_0}+
\[ \mu_{g,1} = \sum_{j\in \N} \sum_{\substack{k \in [j \epsilon, j]}} \sum_{\tau \in \Lambda_{j,k}} \sum_{T \in \mathbb{T}_{\tau,g}} M_T\mu. \]
%and 
%\[ \mu_{g,2} = \sum_{j\in N} \sum_{k\le \epsilon j} 
%        \sum_{\tau \in \Lambda_{j,k}}
%        \mu*\widecheck{\psi_{\tau}}=:\mu_0+\mu_{g, 1}+\mu_{g, 2}. 
%    \end{split}
%\end{equation}
The above used that $\mathbb{T}_{\tau,b}$ is empty when $\tau \in \Lambda_{j,k}$, $k \leq j \epsilon$ and $j$ is sufficiently large, which follows from the Frostman condition on $\mu$.

\begin{claim}\label{220622claim5_6}
Let $\tau\in \Lambda_{j, k}$ with $k<j$. If there exist $\theta\in [0, a]$ and $(\eta_1, \eta_2)$ satisfying 
\begin{equation}
    \eta_1 \gamma'(\theta)+ \eta_2 \left( \gamma \times \gamma'\right)(\theta)\in \supp(\psi_{\tau}), 
\end{equation}
then it holds that $|\eta_2|\sim 2^j$ and $|\eta_1|\sim 2^{j-k/2}$. 
\end{claim}
\begin{proof}[Proof of Claim \ref{220622claim5_6}]
Recalling \eqref{taudefn}, we may assume \begin{multline*}\tau=\tau(\theta',j,k) \\
= \Big\{\lambda_1 \gamma'(\theta') +\lambda_2 \left(\gamma \times \gamma'\right)(\theta')+ \lambda_3 \gamma(\theta') :
 |\lambda_1| \lesssim 2^{-k/2 + j},\lambda_2 \sim 2^{j},  |\lambda_3| \sim 2^{-k + j} \Big\}. \end{multline*}
If 
$\eta_1 \gamma'(\theta)+ \eta_2 \left( \gamma \times \gamma'\right)(\theta)\in \tau(\theta',j,k),$
we discuss some geometric observations.
Noting that $|\lambda_3|\sim 2^{-k+j}$ in the definition of $\tau(\theta',j,k)$, we see that $\eta_1 \gamma'(\theta)+ \eta_2 \left( \gamma \times \gamma'\right)(\theta)\notin \tau(\theta',j,k)$ if $\theta'=\theta$; we also note that if $|\theta'-\theta|\gg 2^{-k/2}$ are too far apart, then $\eta_1 \gamma'(\theta)+ \eta_2 \left( \gamma \times \gamma'\right)(\theta)\notin \tau(\theta',j,k)$. Therefore we must have $|\theta'-\theta|\sim 2^{-k/2}$. In this case, in order for $\eta_1 \gamma'(\theta)+ \eta_2 \left( \gamma \times \gamma'\right)(\theta)\in \tau(\theta',j,k)$, we must have $|\eta_2|\sim 2^j$ and $|\eta_1|\sim 2^{j-k/2}$, which finishes the proof.
\end{proof}
By Claim \ref{220622claim5_6},% and noting $k\le \e j$ in the definition of $\mu_{g,2}$, 
we see that %the contribution from $\mu_{g, 2}$ to 
\eqref{localise} is bounded by 
\begin{align}\label{energymu}
    &1+ \int \int_{\{ |\eta_1| \geq |\eta_2|^{1-\epsilon} \} } \left\lvert \widehat{\mu}\left( \eta_1 \gamma'(\theta)+ \eta_2 \left( \gamma \times \gamma'\right)(\theta) \right)\right\rvert^2 \, d\eta \, d\theta \\
    \label{mug1} &\quad + \int \int_{\{ |\eta_1| < |\eta_2|^{1-\epsilon} \}  } \left\lvert \widehat{\mu_{g,1}}\left( \eta_1 \gamma'(\theta)+ \eta_2 \left( \gamma \times \gamma'\right)(\theta) \right)\right\rvert^2 \, d\eta \, d\theta. \end{align}
For the first term, the change of variables
\begin{equation}\label{220623e5_29}
    \xi = \xi(\eta, \theta) =  \eta_1 \gamma'(\theta)+ \eta_2 \left( \gamma \times \gamma'\right)(\theta)
\end{equation}
has Jacobian 
\begin{align}  \label{stp1} \left\lvert \frac{ \partial(\xi_1, \xi_2, \xi_3) }{\partial(\eta_1, \eta_2, \theta) }(\eta_1, \eta_2, \theta) \right\rvert &=  \left\lvert \eta_1\right\rvert \left\lvert \det\left(\left(\gamma\times \gamma'\right)(\theta), \gamma'(\theta),  \gamma''(\theta) \right) \right\rvert \\
\notag &= |\eta_1| \left\lvert \left\langle \gamma(\theta), \gamma''(\theta) \right\rangle \right\rvert= |\eta_1|, \end{align} 
where in the last step we used 
\begin{equation}
    \inn{\gamma(t)}{\gamma'(t)}\equiv 0\implies \inn{\gamma(t)}{\gamma''(t)}\equiv -1. 
\end{equation}
Applying this change of variables to \eqref{energymu} gives
\begin{align*} \eqref{energymu}\lesssim   1 +  \int_{|\xi|\ge 1} \left\lvert \xi \right\rvert^{\epsilon -1} \left\lvert \widehat{\mu}(\xi)\right\rvert^2 \, d\xi  %&\lesssim   1 +  \int_{|\xi|\ge 1} \left\lvert \xi \right\rvert^{\epsilon -1} \left\lvert \widehat{\mu}(\xi)\right\rvert^2 \, d\xi 
\lesssim  1 +
I_{\alpha-\epsilon}(\mu) \lesssim 1. \end{align*} 
Here $I_{\alpha-\epsilon}(\mu)=\int |\xi|^{\alpha-\e-3}|\wh\mu(\xi)|^2 \, d\xi$ is the $(\alpha-\epsilon)$-energy of $\mu$ and we used the fact that $\alpha>2$ and $\epsilon$ is sufficiently small. The last step is because $c_\alpha(\mu)< \infty$. \\

 It remains to bound the contribution from $\mu_{g, 1}$, in \eqref{mug1}. By frequency disjointness, 
 \begin{equation}\label{220623e5_32}
     \begin{split}
         & \int \int_{\mathbb{R}^2} \left\lvert \widehat{\mu_{g, 1}}\left( \eta_1 \gamma'(\theta)+ \eta_2 \left( \gamma \times \gamma'\right)(\theta) \right)\right\rvert^2 \, d\eta \, d\theta\\
         & \lesssim \sum_{j}\sum_{k\in [\epsilon j, j]} \sum_{\tau\in \Lambda_{j, k}} 
         \int \int_{\R^2} \anorm{
         \sum_{T\in \T_{\tau, g}}
         \widehat{M_T \mu}\left( \eta_1 \gamma'(\theta)+ \eta_2 \left( \gamma \times \gamma'\right)(\theta) \right)
         }^2  
         d\eta \, d\theta. 
     \end{split}
 \end{equation}
 Consider the case $k<j$ and $k=j$ separately. In the former case, we apply the change of variables as in \eqref{220623e5_29} and obtain 
 \begin{multline}\label{220623e5_33}
    \sum_{\tau\in \Lambda_{j, k}}
         \int \int_{\R^2} \anorm{
         \sum_{T\in \T_{\tau, g}}
         \widehat{M_T \mu}\left( \eta_1 \gamma'(\theta)+ \eta_2 \left( \gamma \times \gamma'\right)(\theta) \right)
         }^2  
         d\eta \,  d\theta\\
         \\ 
         \lesssim \sum_{\tau\in \Lambda_{j, k}}
     2^{-j+k/2}
     \int_{\R^3} 
     \anorm{
         \sum_{T\in \T_{\tau, g}}
         \widehat{M_T \mu}(\xi)
         }^2  
         d\xi
         \\ 
         \lesssim \sum_{\tau\in \Lambda_{j, k}}
     \sum_{T\in \T_{\tau, g}}
     2^{-j+k/2} 2^{O(\delta) k}
     \int_{\R^3} 
     \anorm{
         M_T \mu(x)
         }^2  
         dx. \end{multline}
When $k=j$, we show that \eqref{220623e5_33} holds as well. To see this, we first observe that for each fixed $T$, in order for 
\begin{equation}
    \widehat{M_T \mu}\left( \eta_1 \gamma'(\theta)+ \eta_2 \left( \gamma \times \gamma'\right)(\theta) \right)
\end{equation}
not to vanish, $\theta$ has to take values on an interval of length $2^{-j/2}$; next, we apply the two dimensional Plancherel's theorem in the $\eta_1$ and $\eta_2$ variables for every fixed $\theta$, and \eqref{220623e5_33} follows from the uncertainty principle. \\

We continue to estimate \eqref{220623e5_33} and do not distinguish $k<j$ and $k=j$ anymore. We have 
\begin{align}
    \eqref{220623e5_33} & \lesssim \sum_{\tau\in \Lambda_{j, k}}
     \sum_{T\in \T_{\tau, g}}
     2^{-j+k/2} 2^{O(\delta) k}
     \int_{\R^3} 
     \anorm{
         M_T \mu(x)
         }^2  
         dx \label{220623e5_35}\\
    & = 
    2^{-j+k/2} 2^{O(\delta) k}
    \int 
    \sum_{\tau\in \Lambda_{j, k}} 
    \sum_{T\in \T_{\tau, g}} 
    f_T
    d\mu\label{220623e5_36}
\end{align}
where 
\begin{equation}
    f_T:=
    (\eta_T M_T \mu)* \widecheck{\psi_{\tau}},
\end{equation}
and from \eqref{220623e5_35} to \eqref{220623e5_36} we applied Fubini and expanded the square. We cut the unit ball into small balls $B_{\iota}$ of radius $2^{-j+k}$ and let $\nu_{\iota}$ be the restriction of $\mu$ to $2^{10\delta k}B_{\iota}$. By Cauchy-Schwarz, 
\begin{equation}\label{220624e5_38}
    \eqref{220623e5_36} \lesssim 
    2^{-j+k/2} 2^{O(\delta) k}
    \sum_{\iota}
    \mu(2^{10k\delta}B_{\iota})^{1/2} 
    \pnorm{
    \int
    \anorm{
    \sum_{\tau\in \Lambda_{j, k}}
    \sum_{T\in \T_{\tau, g}} 
    f_T
    }^2 d\nu_{\iota}
    }^{1/2}
\end{equation}
Let $\zeta_j$ be a non-negative bump function such that $\widehat{\zeta_j}(\xi)=1$ for $|\xi|\le 2^{j+10}$. By the Fourier support information of $f_T$, we have
\begin{equation}\label{220623e5_39}
    \int
    \anorm{
    \sum_{\tau\in \Lambda_{j, k}}
    \sum_{T\in \T_{\tau, g}} 
    f_T
    }^2 d\nu_{\iota}
    =
    \int
    \anorm{
    \sum_{\tau\in \Lambda_{j, k}}
    \sum_{T\in \T_{\tau, g}} 
    f_T
    }^2 
    d(\nu_{\iota}*\zeta_j )
\end{equation}
By pigeonholing, we can find a subset 
\begin{equation}
    \W_{\iota}\subset \bigcup_{\tau\in \Lambda_{j, k}} \{T\in \T_{\tau, g}: T\cap B_{\iota}\neq \emptyset\}
\end{equation}
 such that $\norm{f_T}_2$ is constant up to a factor of $2$ as $T$ varies over $\W_{\iota}$, and 
 \begin{align}
     \eqref{220623e5_39} \lesssim 
     2^{O(\delta) j}
     \int 
     \anorm{
     \sum_{T\in \W_{\iota}} f_T
     }^2 d(\nu_{\iota}* \zeta_j)
 \end{align}
 By pigeonholing again and by H\"older's inequality, there is a disjoint union $Y$ of balls $Q$ of radius $2^{-j}$, such that 
\begin{equation} \label{holder} 
\int 
\anorm{
\sum_{T \in \mathbb{W}_{\iota}} f_T 
}^2 
d(\nu_{\iota} \ast \zeta_j) \lesssim 
2^{O(\delta) j} 
\norm{
\sum_{T\in \W_{\iota}} f_T
}_{L^p(Y)}^2 
\left( \int_Y \left( \nu_{\iota} \ast \zeta_j\right)^{\frac{p}{p-2}} \right)^{1-\frac{2}{p}},  \end{equation}
and such that each $Q \subseteq Y$ intersects $\sim M$ planks $3T$ as $T$ varies over $\mathbb{W}_{\iota}$, for some dyadic number $M$. By rescaling and then applying the refined decoupling inequality in Theorem~\ref{refineddecoupling} from the Appendix A, the first term in \eqref{holder} satisfies
\begin{align}\label{220624e5_43}
\Norm{
\sum_{T\in \W_{\iota}}
f_T
}_{L^p(Y)}
\lesssim
2^{(3j-\frac{3k}{2})(\frac{1}{2}-\frac{1}{p})+O(\epsilon) j} 
\left( \frac{M}{\left\lvert \mathbb{W}_{\iota} \right\rvert } \right)^{\frac{1}{2}-\frac{1}{p}} 
\left( \sum_{T \in \mathbb{W}_{\iota} } \left\lVert f_T \right\rVert_2^2 \right)^{1/2}. 
\end{align}
For the second term in \eqref{holder}, the assumption that $c_{\alpha}(\mu)$ is finite implies that 
\begin{equation}
    \|\nu_{\iota} \ast \zeta_j\|_{\infty} \lesssim 2^{j(3-\alpha)}.
\end{equation}
Hence by H\"older's inequality and the definition of $\T_{\tau, g}$, we have 
\begin{align} 
\int_Y 
\left( \nu_{\iota} \ast \zeta_j\right)^{\frac{p}{p-2}} 
&\lesssim
2^{\frac{2j(3-\alpha)}{p-2}} M^{-1} \sum_{T \in \mathbb{W}_{\iota}} \left( \nu_{\iota} \ast \zeta_j \right)(3T) \\
&\lesssim 2^{\frac{2j(3-\alpha)}{p-2}} M^{-1} \left\lvert  \mathbb{W}_{\iota} \right\rvert 2^{-k(\alpha_0+1)/2 -\alpha(j-k)}. \label{220624e5_45} 
\end{align} 
Combining \eqref{220624e5_43}  and \eqref{220624e5_45}, we obtain 
\begin{equation} \label{combined} \eqref{holder} 
\lesssim 
2^{O(\epsilon) j}  2^{j(3-\alpha) +k\left(\frac{1}{2}-\frac{1}{p}\right)\left(-4 + 2\alpha - \alpha_0\right)} \sum_{T \in \mathbb{W}_{\iota}} \left\lVert f_T \right\rVert_2^2. \end{equation}
Note that 
\begin{align}
\|f_T\|_2 \lesssim \|M_T \mu \|_2
\end{align}
for every $T$. Substituting into \eqref{combined} and then into \eqref{220623e5_35} yields
\begin{align}
    & \sum_{\tau\in \Lambda_{j, k}} 
    \sum_{T\in \T_{\tau, g}}
    \int_{\R^3}
    \anorm{
    M_T \mu(x)
    }^2
    dx\\
    & \lesssim 
    \sum_{\iota} \mu(2^{10k\delta} B_{\iota})^{1/2}
    2^{O(\epsilon) j} 
    2^{\frac{1}{2} j(3-\alpha) +\frac{1}{2} k\left(\frac{1}{2}-\frac{1}{p}\right)\left(-4 + 2\alpha - \alpha_0\right)} 
    \pnorm{
    \sum_{T \in \mathbb{W}_{\iota}} \left\lVert f_T \right\rVert_2^2
    }^{1/2}.
\end{align}
% \begin{align} 
% \sum_{|j'-j| \leq C } \sum_{\substack{|k'-k| \leq C \\ k ' \leq j'}} \sum_{\tau \in \Lambda_{j',k'}} \sum_{\substack{T \in \mathbb{T}_{\tau, g}  \\ T \cap B_m \neq \emptyset }} \int \left\lvert M_T \mu \right\rvert^2 \\
% \lesssim \mu\left( 2^{100k\delta} B_m \right)^{1/2} 2^{10j\epsilon+ \frac{1}{2} \left[ j(3-\alpha) + k\left(\frac{1}{2}-\frac{1}{p}\right)\left(-4 + 2\alpha - \alpha_0 \right) \right] } \\
% \times \left( \sum_{|j'-j| \leq C } \sum_{\substack{|k'-k| \leq C \\ k ' \leq j'}} \sum_{\tau \in \Lambda_{j',k'}} \sum_{\substack{T \in \mathbb{T}_{\tau, g}  \\ T \cap B_m \neq \emptyset }} \int \left\lvert M_T \mu \right\rvert^2  \right)^{1/2}. 
% \end{align}  
By Cauchy-Schwarz in the sum over $\iota$, we obtain 
\begin{align}
    & \sum_{\tau\in \Lambda_{j, k}} 
    \sum_{T\in \T_{\tau, g}}
    \int_{\R^3}
    \anorm{
    M_T \mu(x)
    }^2
    dx
    \lesssim 
    2^{O(\epsilon) j} 
    2^{j(3-\alpha) + k\left(\frac{1}{2}-\frac{1}{p}\right)\left(-4 + 2\alpha - \alpha_0\right)}.
\end{align}
By substituting back into \eqref{220623e5_35}, we obtain 
\begin{equation}
    \eqref{220623e5_33} 
    \lesssim 
    \sum_j \sum_{\epsilon j\le k\le j} 2^{-j+k/2} 2^{O(\epsilon) j} 
    2^{j(3-\alpha) + k\left(\frac{1}{2}-\frac{1}{p}\right)\left(-4 + 2\alpha - \alpha_0\right)}.
\end{equation}
In the end, we pick $p=4$ and finish the proof. \end{proof}

\appendix 
 
\section{Refined decoupling inequality}

The refined decoupling inequality stated here is a natural analogue of the refined decoupling inequality for the paraboloid from \cite{GIOW}. The shortest length of a plank dual to an $R^{-1/2}$-cap in the cone is $\approx 1$, rather than $\approx R^{1/2}$ in the case of the paraboloid, so the setup uses unit cubes instead of $R^{1/2}$-cubes. The argument is similar to the paraboloid case, using induction and with Lorentz rescaling in place of parabolic rescaling, but the use of unit cubes requires the induction to be carried out over a finer sequence of scales (similarly to the induction setup in \cite{duzhang}). The refined decoupling inequalities in \cite{harris2019improved, harris2021restricted} used tubes rather than planks, and the use of planks here is a significant reason for the improved result on the projection problem (at least in Theorem~\ref{thmpositive}). Much of the proof is similar to \cite{harris2019improved, harris2021restricted}, with only the differences outlined above.

For each $R \geq 1$ let $\Xi_R = \left\{ j R^{-1/2}  : j \in \mathbb{Z} \right\} \cap [0, \aaa].$
For each $\theta \in \Xi_R$, let 
\begin{equation} 
\tau_R(\theta) =\bigg\{ x_1 \gamma(\theta) + x_2 \gamma'(\theta)+ x_3 (\gamma \times \gamma' )(\theta ):
 1 \leq x_1 \leq 2, \, |x_2| \leq R^{-1/2}, \, |x_3| \leq R^{-1}  \bigg\}. 
 \end{equation}
If it is clear from the context which $R$ is used, then we often abbreviate $\tau_R(\theta)$ to $\tau(\theta)$. Let $\mathcal{P}_{R^{-1}}= \left\{ \tau(\theta): \theta\in \Xi_R \right\}.$ For $\tau=\tau(\theta)\in \mathcal{P}_{R^{-1}}$, denote $\theta_{\tau}:=\theta$. Let 
\begin{equation} 
T^{\circ}_{\tau,0} = \bigg\{ x_1\left(\gamma \times \gamma' \right)(\theta_{\tau} )+ x_2 \gamma'(\theta_{\tau})+ x_3 \gamma(\theta_{\tau} ) :
|x_1| \leq R, \, |x_2| \leq R^{1/2}, \, |x_3| \leq 1  \bigg\}. 
\end{equation}
Moreover, we will use $\T^{\circ}_{\tau}$ to denote the collection of translates of $T^{\circ}_{\tau, 0}$ that cover $B(0, R)$. 
For a fixed small constant $\delta>0$, denote $T_{\tau, 0}:=R^{\delta} T^{\circ}_{\tau, 0}$, and $\T_{\tau}:=\{R^{\delta} T: T\in \T^{\circ}_{\tau}\}$. For $T\in \T_{\tau}$, set $\tau(T)=\tau$. 
\begin{definition}
Fix $T\in \T_{\tau}$. We say that a function $f_T: \R^3\to \C$ is a $T$-function if $\widehat{f_T}$ is supported on $\tau(T)$ and 
\begin{equation}\label{microlocalised}  \left\lVert f_T \right\rVert_{L^{\infty} \left( B(0,R) \setminus T \right) }\lesssim_{\delta} R^{-10000} \left\lVert f_T \right\rVert_2.
\end{equation}
\end{definition}

 \begin{theorem} \label{refineddecoupling} 
%  Let $\gamma: [a,b] \to\ZS^2$ be a $C^3$ curve with $\det(\gamma, \gamma', \gamma'')$ nonvanishing. Let $B \geq 1$ be such that 
% \begin{equation} \label{gammacdn1} \left\lvert \det\left( \gamma, \gamma', \gamma'' \right) \right\rvert \geq B^{-1} , \end{equation}
% and 
% \begin{equation} \label{gammacdn2} \left\lVert \gamma\right \rVert_{C^3[a,b]} \leq B. \end{equation}
Let $\gamma: [a,b] \to \ZS^2$ be a $C^2$ curve with $\det(\gamma, \gamma', \gamma'')$ nonvanishing. Let $B \geq 1$ be such that 
\begin{equation} \label{gammacdn1} \left\lvert \det\left( \gamma, \gamma', \gamma'' \right) \right\rvert \geq B^{-1} , \end{equation}
and 
\begin{equation} \label{gammacdn2} \left\lVert \gamma\right \rVert_{C^2[a,b]} \leq B. \end{equation}
Let $R \geq 1$ and suppose that
\[ f = \sum_{T \in \mathbb{W}} f_T,  \]
where each $f_T$ is a $T$-function and 
\[ \mathbb{W} \subseteq  \bigcup_{\tau \in \mathcal{P}_{R^{-1}}}  \mathbb{T}_{\tau}. \]
Assume that for all $T, T' \in \mathbb{W}$, 
\begin{equation}
    \norm{f_T}_2\sim \norm{f_{T'}}_2.
\end{equation}
Let $Y$ be a disjoint union of unit balls in $B(0,R)$, each of which intersects at most $M$ sets $2T$ with $T \in \mathbb{W}$. Then for $2 \leq p \leq 6$, 
\[ \left\lVert f \right\rVert_{L^p(Y)} \lesssim_{\epsilon, \delta} R^{\epsilon} \left( \frac{M}{\left\lvert \mathbb{W} \right\rvert } \right)^{\frac{1}{2} - \frac{1}{p} } \left( \sum_{T \in \mathbb{W}} \left\lVert f_T \right\rVert_p^2 \right)^{1/2}. \] \end{theorem}
\begin{proof} Assume that $[a,b] = [-1,1]$. Fix $\epsilon \in (0,1/2)$, $\delta_0 = \epsilon^{100}$, $\delta \in (0, \delta_0)$, 
\[ R \geq \min\left\{B^{10^3/\epsilon}, 2^{10^5/\epsilon}\right\}, \]
and assume inductively that a (superficially) stronger version of the theorem holds with $K^2$-cubes instead of unit cubes, where $K = R^{\delta^2}$, for all scales smaller than $\widetilde{R} := R/K^2$, for all curves $\gamma$ satisfying \eqref{gammacdn1} and \eqref{gammacdn2}, and for all $B \geq 1$. 

For each $\tau \in \mathcal{P}_{R^{-1/2}}(\Gamma(\gamma))$, let $\kappa = \kappa(\tau) \in \mathcal{P}_{K^{-1} }(\Gamma(\gamma))$ be the element of $\mathcal{P}_{K^{-1} }(\Gamma(\gamma))$ which minimises $|\theta_{\tau} - \theta_{\kappa} |$. For each $\kappa$, let 
\begin{multline*} \Box_{\kappa, 0} = \Big\{ x_1 \frac{\left(\gamma \times \gamma' \right)(\theta_{\kappa} )}{\left\lvert \left(\gamma \times \gamma' \right)(\theta_{\kappa} ) \right\rvert }+ x_2 \frac{\gamma'(\theta_{\kappa} )}{\left\lvert \gamma'(\theta_{\kappa} ) \right\rvert} + x_3 \gamma(\theta_{\kappa} ) : \\
|x_1| \leq R^{1+\delta}, \, |x_2| \leq R^{1+\delta}/K, \, |x_3| \leq R^{1+\delta}/K^2  \Big\}, \end{multline*}
and
\begin{multline*} \mathbb{P}_{\kappa} = \Big\{ \Box = a \gamma(\theta_{\kappa} )+  b \frac{\gamma'(\theta_{\kappa} )}{\left\lvert \gamma'(\theta_{\kappa} ) \right\rvert} +\Box_{\kappa, 0} :  \\
  a\in \left((1/10) R^{1+\delta}K^{-2}\right) \mathbb{Z}, \quad b \in \left( (1/10) R^{1+\delta}K^{-1} \right) \mathbb{Z} \Big\}.  \end{multline*} 
Let $\mathbb{P} = \bigcup_{\kappa \in \mathcal{P}_{K^{-1} }(\Gamma(\gamma))} \mathbb{P}_{\kappa}$. Given any $\tau$ and corresponding $\kappa = \kappa(\tau)$, 
\begin{equation} \label{angledistortion} \left\lvert \left\langle (\gamma \times \gamma')(\theta_{\tau}), \gamma'(\theta_{\kappa}) \right\rangle \right\rvert \leq B^{-7} K^{-1}, \end{equation}
and 
\begin{equation} \label{angledistortion2} \left\lvert \left\langle (\gamma \times \gamma')(\theta_{\tau}), \gamma(\theta_{\kappa}) \right\rangle \right\rvert \leq B^{-7}K^{-2}. \end{equation}
 It follows that for each $T \in \mathbb{T}_{\tau}$, there are $\sim 1$ sets $\Box \in \mathbb{P}_{\kappa(\tau)}$ with $T \cap 10\Box \neq \emptyset$, and moreover $T \subseteq 100\Box$ whenever $T \cap 10\Box \neq \emptyset$. For each such $T$ let $\Box = \Box(T)\in \mathbb{P}_{\kappa}$ be some choice such that $T \cap 10\Box \neq \emptyset$, and let $\mathbb{W}_{\Box}$ be the set of $T$'s associated to $\Box$.

For each $\kappa$ and $\Box \in \mathbb{P}_{\kappa}$, let $\left\{Q_{\Box} \right\}_{Q_{\Box}}$ be a finitely overlapping cover of $100\Box$ by translates of the ellipsoid
\begin{multline*} \Bigg\{ x_1 \gamma(\theta_{\kappa} ) + x_2 \frac{\gamma'(\theta_{\kappa} )}{\left\lvert \gamma'(\theta_{\kappa} )\right\rvert} + x_3 \frac{\left( \gamma \times \gamma' \right)(\theta_{\kappa} )}{\left\lvert \left( \gamma \times \gamma' \right)(\theta_{\kappa} )\right\rvert}: \\
 \left( |x_1| ^2 + \left(|x_2| K^{-1}\right)^2 + \left(|x_3 | K^{-2} \right)^2 \right)^{1/2} \leq  \widetilde{K}^2 \Bigg\}. \end{multline*} 
Using Poisson summation, let $\{ \eta_{Q_{\Box}} \}_{Q_{\Box} \in \mathcal{Q}_{\Box}}$ be a smooth partition of unity such that on $10^3 \Box$,
\[ \sum_{Q_{\Box} \in \mathcal{Q}_{\Box}} \eta_{Q_{\Box}} =1, \]
and such that each $\eta_{Q_{\Box}}$ satisfies 
\[ \|\eta_{Q_{\Box}}\|_{\infty} \lesssim 1, \quad  \|\eta_{Q_{\Box}}\|_{L^{\infty}(\mathbb{R}^3 \setminus Q_{\Box} ) } \lesssim R^{-10000}, \]
and 
\[  |\eta_{Q_{\Box} }(x) | \lesssim  \dist(x, Q_{\Box})^{-10000} \quad \forall \, x \in \mathbb{R}^3,  \]
with $\widehat{\eta_{Q_{\Box}}}$ supported in 
\begin{multline*} \bigg\{ \xi_1 \gamma(\theta_{\kappa} ) + \xi_2 \frac{\gamma'(\theta_{\kappa} )}{\left\lvert \gamma'(\theta_{\kappa} )\right\rvert} + \xi_3 \frac{\left( \gamma \times \gamma' \right)(\theta_{\kappa} )}{\left\lvert \left( \gamma \times \gamma' \right)(\theta_{\kappa} )\right\rvert} \\
:  |\xi_1| \leq \widetilde{K}, \quad |\xi_2| \leq \widetilde{K} K, \quad |\xi_3| \leq \widetilde{K} K^2 \bigg\}. \end{multline*}
  By dyadic pigeonholing,  
\[ \left\lVert f \right\rVert_{L^p(Y)} \lesssim \log R  \left\lVert \sum_{\Box} \sum_{T \in \mathbb{W}_{\Box}} \eta_{Y_{\Box}} f_T \right\rVert_{L^p(Y)} + R^{-1000} \left( \sum_{T \in \mathbb{W} } \left\lVert f_T \right\rVert_p^2 \right)^{1/2}, \]
where, for each $\Box$, $Y_{\Box}$ is a union over a subset of the sets $Q_{\Box}$, and $\eta_{Y_{\Box}}$ is the corresponding sum over $\eta_{Q_{\Box}}$, such that each $Q_{\Box} \subseteq Y_{\Box}$ intersects a number $\# \in [M'(\Box), 2M'(\Box))$ different sets $(1.5)T$ with $T \in \mathbb{W}_{\Box}$, up to a factor of 2. By pigeonholing again,
\[ \left\lVert \sum_{\Box} \sum_{T \in \mathbb{W}_{\Box}} \eta_{Y_{\Box}} f_T \right\rVert_{L^p(Y)} \lesssim (\log R )^2 \left\lVert \sum_{\Box \in \mathbb{B}} \sum_{T \in \mathbb{W}_{\Box}} \eta_{Y_{\Box}} f_T \right\rVert_{L^p(Y)}, \]
where $\left\lvert \mathbb{W}_{\Box} \right\rvert$ and $M'=M'(\Box)$ are constant over $\Box \in \mathbb{B}$, up to a factor of 2. By one final pigeonholing step,
\[ \left\lVert \sum_{\Box \in \mathbb{B}} \sum_{T \in \mathbb{W}_{\Box}} \eta_{Y_{\Box}} f_T\right\rVert_{L^p(Y)} \lesssim \log R \left\lVert \sum_{\Box \in \mathbb{B}} \sum_{T \in \mathbb{W}_{\Box}} \eta_{Y_{\Box}} f_T \right\rVert_{L^p(Y')}, \] 
where $Y'$ is a union over $K^2$-balls $Q \subseteq Y$ such that each ball $2Q$ intersects a number $\# \in [M'', 2M'')$ of the sets $Y_{\Box}$ in a set of strictly positive Lebesgue measure, as $\Box$ varies over $\mathbb{B}$. Fix $Q \subseteq Y'$. By the decoupling theorem for generalised $C^2$ cones, followed by Hölder's inequality,
\begin{multline*} \left\lVert \sum_{\Box \in \mathbb{B}} \sum_{T \in \mathbb{W}_{\Box}} \eta_{Y_{\Box}} f_T \right\rVert_{L^p(Q)} \\
\leq C_{\epsilon} B^{100} K^{\epsilon/100} \left( M'' \right)^{\frac{1}{2} - \frac{1}{p}} \left(  \sum_{\Box \in \mathbb{B}} \left\lVert\sum_{T \in \mathbb{W}_{\Box}} \eta_{Y_{\Box}} f_T \right\rVert_{L^p(2Q)}^p \right)^{1/p} \\
+ R^{-900} \left( \sum_{T \in \mathbb{W}} \left\lVert f_T \right\rVert_p^2 \right)^{1/2}. \end{multline*}
Summing over $Q$ gives 
\begin{multline*} \|f\|_{L^p(Y)} \lesssim C_{\epsilon}\left( \log R \right)^{100} B^{100} K^{\epsilon/100} \left( M'' \right)^{\frac{1}{2} - \frac{1}{p}} \\
\times \left(  \sum_{\Box \in \mathbb{B}} \left\lVert\sum_{T \in \mathbb{W}_{\Box}} f_T \right\rVert_{L^p(Y_{\Box} ) }^p \right)^{1/p} + R^{-800} \left( \sum_{T \in \mathbb{W}} \left\lVert f_T \right\rVert_p^2 \right)^{1/2}. \end{multline*} This will be bounded using the inductive assumption, following a Lorentz rescaling.

For each $\theta \in [-1,1]$, define the Lorentz rescaling map $L = L_{\theta}$ at $\theta$ by
\begin{multline*} L\left[x_1 \gamma(\theta) + x_2 \frac{\gamma'(\theta)}{\left\lvert \gamma'(\theta)\right\rvert}  + x_3 \frac{\left(\gamma \times \gamma' \right)(\theta)}{\left\lvert \left(\gamma \times \gamma' \right)(\theta)\right\rvert}\right] \\
=x_1 \gamma(\theta) + K x_2  \frac{\gamma'(\theta)}{\left\lvert \gamma'(\theta)\right\rvert}  + K^2 x_3 \frac{\left(\gamma \times \gamma' \right)(\theta)}{\left\lvert \left(\gamma \times \gamma' \right)(\theta)\right\rvert}. \end{multline*}
Let 
\[ \widetilde{\gamma}(\phi) = \frac{ L(\gamma(\phi))}{\left\lvert L(\gamma(\phi) ) \right\rvert}, \qquad \phi \in [-1,1]. \]
Then for any $\phi \in [-1,1]$,
\[ \widetilde{\gamma}'(\phi) = \frac{\pi_{\widetilde{\gamma}(\phi)^{\perp} }\left( L(\gamma'(\phi) ) \right) }{\left\lvert L(\gamma(\phi) ) \right\rvert}, \]
and 
\[ \widetilde{\gamma}''(\phi) = \frac{\pi_{\widetilde{\gamma}(\phi)^{\perp} }\left( L(\gamma''(\phi) ) \right) }{\left\lvert L(\gamma(\phi) ) \right\rvert} - \frac{\left\langle L(\gamma(\phi)), L(\gamma'(\phi) ) \right\rangle \pi_{\widetilde{\gamma}(\phi)^{\perp} }\left( L(\gamma'(\phi) ) \right) }{\left\lvert L(\gamma(\phi) ) \right\rvert^3}. \]
Hence 
\begin{align*} \det( \widetilde{\gamma}, \widetilde{\gamma}',\widetilde{\gamma}'' ) &= \frac{1}{\left\lvert L \circ \gamma \right\rvert^3}  \det\left( L \circ \gamma, L \circ \gamma', L \circ \gamma''\right) \\
&= \frac{K^3}{\left\lvert L \circ \gamma \right\rvert^3}  \det\left( \gamma, \gamma', \gamma''\right). \end{align*}
Let $\varepsilon = (10^5 B^{10})^{-1}$, and for fixed $\theta \in [-1+\varepsilon,1-\varepsilon]$, let
\begin{equation} \label{scalesize} \widehat{\gamma}(\phi) = \widetilde{\gamma}(\theta +   K^{-1}\phi  ), \quad  \phi \in [-\varepsilon, \varepsilon]. \end{equation}
The assumption that $\left\lVert \gamma\right \rVert_{C^2[-1,1]} \leq B$ yields
\[ 1 \leq |L(\gamma(\phi) ) | \leq 1+ 10B \varepsilon, \quad \forall \, \phi \in [\theta - \varepsilon K^{-1}, \theta + \varepsilon K^{-1} ]. \]
Similarly, 
\[ \left\lvert L(\gamma'(\phi)) - L(\gamma'(\theta)) \right\rvert \leq 10\varepsilon B K,  \quad \forall \, \phi \in [\theta - \varepsilon K^{-1}, \theta + \varepsilon K^{-1} ]. \]
%\[ \left\lvert L(\gamma''(\phi)) - L(\gamma''(\theta) ) \right\rvert \leq 10\varepsilon B K,  \quad \forall \, \phi \in [\theta - \varepsilon K^{-1}, \theta + \varepsilon K^{-1} ],  \]
%and 
%\[|L \gamma'''(\phi) | \leq BK^2, \quad \forall \, \phi \in [\theta - \varepsilon K^{-1}, \theta + \varepsilon K^{-1} ]. \]
It follows that
\[ \left\lvert \det\left(\widehat{\gamma}, \widehat{\gamma}', \widehat{\gamma}''\right) \right\rvert \geq (2B)^{-1} \]
on $[-\varepsilon,\varepsilon]$, and that
\[ \left\lVert \widehat{\gamma} \right\rVert_{C^2[-\varepsilon,\varepsilon] } \leq 2B. \]

For each $\Box \in \mathbb{B}$, given $T \in \mathbb{W}_{\Box}$, let $g_T =  f_T \circ L$, where $L= L_{\theta_{\kappa(\Box) } }$. Then
\begin{equation} \label{inductthis} \left\lVert\sum_{T \in \mathbb{W}_{\Box}} f_T \right\rVert_{L^p(Y_{\Box} ) } \leq K^{\frac{3}{p}}\left\lVert\sum_{T \in \mathbb{W}_{\Box}} g_T \right\rVert_{L^p(L^{-1} Y_{\Box} ) }. \end{equation}
The inequalities \eqref{angledistortion} and \eqref{angledistortion2} imply that for each $T \in \mathbb{W}_{\Box}$, the set $L^{-1}(T)$ is a equivalent (up to a factor 1.01) to a plank of length $\widetilde{R}^{1+\delta}$ in its longest direction parallel to $L^{-1}(\gamma \times \gamma')(\theta_{\tau(T)})$, of length $\widetilde{R}^{1/2+\delta}$ in its medium direction, and of length $\widetilde{R}^{\delta}$ in its shortest direction. The ellipsoids $Q_{\Box}$ are rescaled to $\widetilde{K}^2$-balls $L^{-1}(Q_{\Box})$. Moreover, it will be shown that
\begin{multline} \label{scaledtau} L(\tau) \subseteq \bigg\{ x_1 \widetilde{\gamma}(\theta_{\tau} ) + x_2 \frac{\widetilde{\gamma}'(\theta_{\tau})}{\left\lvert \widetilde{\gamma}'(\theta_{\tau}) \right\rvert}  + x_3 \frac{ \left( \widetilde{\gamma}  \times \widetilde{\gamma}'\right)(\theta_{\tau})}{\left\lvert \left( \widetilde{\gamma}  \times \widetilde{\gamma}'\right)(\theta_{\tau})\right\rvert } \\
: 1 \leq x_1 \leq 2.01, \, |x_2| \leq (1.01)\widetilde{R}^{-1/2}, \, |x_3| \leq \widetilde{R}^{-1} \bigg\}. \end{multline} 
To prove this, let
\[ x =  x_1 \gamma(\theta_{\tau}) + x_2 \frac{\gamma'(\theta_{\tau})}{\left\lvert\gamma'(\theta_{\tau}) \right\rvert} + x_3 \frac{(\gamma \times \gamma')(\theta_{\tau})}{\left\lvert (\gamma \times \gamma')(\theta_{\tau})\right\rvert} \in \tau, \]
where 
\[ x_1 \in [1,2], \quad |x_2| \leq R^{-1/2}, \quad |x_3| \leq R^{-1}. \]
The vector $\left(\widetilde{\gamma} \times \widetilde{\gamma}'\right)(\theta_{\tau})$ is parallel to $L^{-1}((\gamma \times \gamma')(\theta_{\tau}))$, since $L^{-1}((\gamma \times \gamma')(\theta_{\tau}))$ is orthogonal to $\widetilde{\gamma}(\theta_{\tau})$ and $\widetilde{\gamma}'(\theta_{\tau})$. The inequality
\[ \left\lvert L^{-1}((\gamma \times \gamma')(\theta_{\tau})) \right\rvert \geq K^{-2}\left\lvert  (\gamma \times \gamma')(\theta_{\tau}) \right\rvert \]
gives
\begin{equation} \label{component1} \left\lvert \left\langle Lx, \frac{ L^{-1}((\gamma \times \gamma')(\theta_{\tau}) )}{\left\lvert  L^{-1}((\gamma \times \gamma')(\theta_{\tau}) ) \right\rvert }\right\rangle \right\rvert \leq \widetilde{R}^{-1}. \end{equation}
Moreover, 
\begin{align} \notag &\left\lvert \left\langle Lx, \frac{ \pi_{L(\gamma(\theta_{\tau}))^{\perp}} \left( L\left(\gamma'(\theta_{\tau}) \right)\right)}{\left\lvert \pi_{L(\gamma(\theta_{\tau}))^{\perp}} \left( L\left(\gamma'(\theta_{\tau}) \right)\right) \right\rvert} \right\rangle \right\rvert \\
\notag & \quad = \left\lvert \left\langle x_2 L\left(\frac{\gamma'(\theta_{\tau})}{\left\lvert\gamma'(\theta_{\tau}) \right\rvert}\right) + x_3 L\left(\frac{(\gamma \times \gamma')(\theta_{\tau})}{\left\lvert (\gamma \times \gamma')(\theta_{\tau})\right\rvert}\right), \frac{ \pi_{L(\gamma(\theta_{\tau}))^{\perp}} \left( L\left(\gamma'(\theta_{\tau}) \right)\right)}{\left\lvert \pi_{L(\gamma(\theta_{\tau}))^{\perp}} \left( L\left(\gamma'(\theta_{\tau}) \right)\right) \right\rvert} \right\rangle \right\rvert \\
\label{component2} &\quad \leq (1.01)\widetilde{R}^{-1/2}.
 \end{align} 
For the direction $L(\gamma(\theta_{\tau}))$, 
\begin{equation} \label{component3} \left\langle Lx, \frac{ L(\gamma(\theta_{\tau}) )}{|L(\gamma(\theta_{\tau}) ) | } \right\rangle = x_1|L(\gamma(\theta_{\tau} ) ) | + O(K^2R^{-1/2}). \end{equation}
Combining \eqref{component1}, \eqref{component2} and \eqref{component3} gives \eqref{scaledtau}. 

Inductively applying the theorem at scale $\widetilde{R}$ gives
\begin{align*} \eqref{inductthis} &\lesssim C_{\epsilon, \delta} B^{10^{10}/\epsilon} R^{\epsilon} K^{-2\epsilon} \left( \frac{ M' }{\left\lvert \mathbb{W}_{\Box} \right\rvert } \right)^{\frac{1}{2} - \frac{1}{p} } \left( \sum_{T \in \mathbb{W}_{\Box} } \left\lVert f_T\right\rVert_p^2 \right)^{1/2} \end{align*} 
for each $\Box \in \mathbb{B}$. Hence 
\begin{multline*} \|f\|_{L^p(Y)}  \\
\leq C_{\epsilon, \delta} B^{10^{10}/\epsilon} K^{-\epsilon} \left( \frac{ M'M''}{\left\lvert \mathbb{W}_{\Box} \right\rvert } \right)^{\frac{1}{2} - \frac{1}{p} } \left(\sum_{\Box \in \mathbb{B}} \left( \sum_{T \in \mathbb{W}_{\Box} } \left\lVert f_T \right\rVert_p^2 \right)^{p/2}\right)^{1/p}. \end{multline*} 
By the dyadically constant property of $\|f_T\|_p$, this is 
\[ \lesssim C_{\epsilon, \delta} B^{10^{10}/\epsilon} K^{-\epsilon} \left( \frac{ M'M'' }{\left\lvert \mathbb{W} \right\rvert  } \right)^{\frac{1}{2} - \frac{1}{p} } \left( \frac{ \left\lvert \mathbb{B} \right\rvert \left\lvert \mathbb{W}_{\Box} \right\rvert }{\left\lvert \mathbb{W}\right\rvert  } \right)^{\frac{1}{p} }\left( \sum_{T \in \mathbb{W} } \left\lVert f_T \right\rVert_2^2 \right)^{1/2}. \]
The second bracketed term is $\lesssim 1$, since
\[ \left\lvert \mathbb{W}  \right\rvert  = \sum_{T \in \mathbb{W}} 1  \geq \sum_{\Box \in \mathbb{B}} \sum_{\substack{ T \in \mathbb{W}: \\
\Box = \Box(T)}} 1 \sim \sum_{\Box \in \mathbb{B}} \sum_{\substack{ T \in \mathbb{W}: \\
\Box = \Box(T)}} 1 \geq \left\lvert \mathbb{B} \right\rvert \left\lvert \mathbb{W}_{\Box} \right\rvert. \]
It remains to show that $M'M'' \lesssim M$. Let $Q \subseteq Y'$ be any $R^{1/2}$-ball. By definition of $M$,
\begin{align*}  M &\gtrsim \sum_{\substack{ T \in \mathbb{W}:  \\
2T \cap Q \neq \emptyset }}  \sum_{\substack{\Box \in \mathbb{B}: \\ 
\Box = \Box(T)}} 1  \\
&= \sum_{\Box \in \mathbb{B}} \sum_{\substack{ T\in \mathbb{W}: \\
\Box = \Box(T) \\
2T \cap Q \neq \emptyset}} 1 \\
&\geq \sum_{\Box \in \mathbb{B}} \sum_{\substack{T \in \mathbb{W}_{\Box}: \\
2T \cap Q \neq \emptyset }} 1. \end{align*}  
By definition of $M'$ and $M''$, \begingroup
\allowdisplaybreaks
\begin{align} \notag M'M'' &\sim \sum_{\substack{\Box \in \mathbb{B}: \\ m(Y_{\Box} \cap 2Q)>0 }} M' \\
\notag &\leq  \sum_{\substack{\Box \in \mathbb{B}: \\ m(Y_{\Box} \cap 2Q)>0 }} \sum_{Q_{\Box} \subseteq Y_{\Box}} M' \frac{m(Q_{\Box} \cap 2Q)}{m(Y_{\Box} \cap 2Q) } \\ 
\notag &\sim  \sum_{\substack{\Box \in \mathbb{B}: \\ m(Y_{\Box} \cap 2Q)>0 }} \sum_{Q_{\Box} \subseteq Y_{\Box}} \sum_{\substack{T\in \mathbb{W}_{\Box}: \\
 Q_{\Box} \cap (1.5)T \neq \emptyset }}  \frac{m(Q_{\Box} \cap 2Q)}{m(Y_{\Box} \cap 2Q) } \\
\notag &=  \sum_{\substack{\Box \in \mathbb{B}: \\ m(Y_{\Box} \cap 2Q)>0 }} \sum_{T \in \mathbb{W}_{\Box}} \sum_{\substack{Q_{\Box} \subseteq Y_{\Box}: \\ Q_{\Box} \cap (1.5)T \neq \emptyset} }  \frac{m(Q_{\Box} \cap 2Q)}{m(Y_{\Box} \cap 2Q) } \\
\label{nontrivial} &\leq  \sum_{\substack{\Box \in \mathbb{B}: \\ m(Y_{\Box} \cap 2Q)>0 }} \sum_{ \substack{T\in \mathbb{W}_{\Box}: \\ 2T \cap Q \neq \emptyset }} \sum_{Q_{\Box} \subseteq Y_{\Box}}  \frac{m(Q_{\Box} \cap 2Q)}{m(Y_{\Box} \cap 2Q) } \\
\notag &\lesssim  \sum_{\Box \in \mathbb{B}: } \sum_{ \substack{T \in \mathbb{W}_{\Box}: \\ 2T \cap Q \neq \emptyset  }} 1 \\
\notag &\lesssim M.  \end{align} \endgroup 
The inequality \eqref{nontrivial} above follows from the observation that if $Q_{\Box} \cap 2Q \neq \emptyset$, and if $T \in \mathbb{W}_{\Box}$ is such that $Q_{\Box} \cap (1.5)T \neq \emptyset$, then $2T \cap Q \neq \emptyset$.  \end{proof}

%\appendix

\section{Decoupling for \texorpdfstring{$C^2$}{C2} cones}

We have been using the decoupling inequality for $C^2$ cones in $\R^3$ in a few places above, but it may have not been written down in the literature. In the appendix, we state it and sketch its proof. We start with the decoupling for $C^2$ curves on $\R^2$.

\begin{theorem}\label{thm10}
Let $\gamma: [-1, 1]\to \R$ with $\gamma(0)=\gamma'(0)=0$ be $C^2$ and satisfy $\gamma''(t)\neq 0$ for every $t\in [-1, 1]$. Then 
\begin{equation}
    \Norm{
    E_{[-1, 1]} f
    }_{L^6(\R^2)} \lesim_{\gamma, \epsilon} \delta^{-\epsilon} 
    \pnorm{
    \sum_{I\subset [0, 1], |I|=\delta}
\norm{E_I f}_{L^6(\R^2)}^2
    }^{1/2},
\end{equation}
for every $\epsilon>0$ and $\delta\in (0, 1)$. 
Here 
\begin{equation}
    E_I f(x, y):=\int_I f(t)e^{i(xt+y\gamma(t))}dt,
\end{equation}
for an interval $I\subset [-1, 1]$. 
\end{theorem}

One can follow the same argument as in \cite{guo2021short} to prove Theorem \ref{thm10}. We leave out the proof. Via the bootstrapping argument as in Bourgain and Demeter \cite{bourgain2015proof} and Pramanik and Seeger \cite{MR2288738}, one can prove the following decoupling for $C^2$ cones. 

\begin{theorem}\label{thm11}
Let $\gamma: [-1, 1]\to \R$ with $\gamma(0)=\gamma'(0)=0$ be $C^2$ and satisfy $\gamma''(t)\neq 0$ for every $t\in [-1, 1]$. Then 
\begin{equation}
    \Norm{
    \mc{E}_{[-1, 1]} f
    }_{L^6(\R^3)} \lesim_{\gamma, \epsilon} \delta^{-\epsilon} 
    \pnorm{
    \sum_{I\subset [0, 1], |I|=\delta}
\norm{\mc{E}_I f}_{L^6(\R^3)}^2
    }^{1/2},
\end{equation}
for every $\epsilon>0$ and $\delta\in (0, 1)$. Here 
\begin{equation}
    \mc{E}_I f(x, y, z):=\int_{[1, 2]\times I} f(s, t)e^{i(sx+st y+s\gamma(t)z)}dsdt,
\end{equation}
for an interval $I\subset [-1, 1]$. 
\end{theorem}

We give a sketch of the proof of Theorem \ref{thm11}. By the triangle inequality, we can assume that $f$ is supported on $[1, 1+\delta^{\epsilon}]\times [0, \delta^{\epsilon}]$. The key observation in \cite{bourgain2015proof} is that the cone 
\begin{equation}\label{cone1z}
    \{s(1, t, \gamma(t)): 1\le s\le 1+\delta^{\epsilon}, 0\le t\le \delta^{\mu}\}
\end{equation}
is in the $\delta^{2\mu+\epsilon}$-neighborhood of the cylinder 
\begin{equation}\label{cylinder1z}
    \{(s, t, \gamma(t)): 1\le s\le 1+\delta^{\epsilon}, 0\le t\le \delta^{\mu}\},
\end{equation}
for every $\mu\ge \epsilon$. To see this, let us take one point $s(1, t, \gamma(t))$ from \eqref{cone1z}, and we will show that its distance to $(s, st, \gamma(st))$, which lies in \eqref{cylinder1z}, is $\lesim \delta^{2\mu+\epsilon}$. This amounts to proving 
\begin{equation}\label{e154}
    |\gamma(st)-s\gamma(t)|\lesim \delta^{2\mu+\epsilon}. 
\end{equation}
Note that 
\begin{equation}
    |\gamma(st)-s\gamma(t)| \lesim |1-s||\gamma(t)|+|\gamma(st)-\gamma(t)|.
\end{equation}
The desired bound \eqref{e154} follows from Taylor's expansion and mean value theorems. 

After proving \eqref{e154}, one can then apply Theorem \ref{thm10} iteratively, in the same way as in \cite{bourgain2015proof}, and finish the proof of Theorem \ref{thm11}. We leave out the iteration step.

\bibliographystyle{abbrv}
\bibliography{bibli}

\begin{thebibliography}{10}

\bibitem{bourgain2015proof}
J.~Bourgain and C.~Demeter.
\newblock The proof of the $\ell^2$ decoupling conjecture.
\newblock {\em Ann. of Math.}, pages 351--389, 2015.

\bibitem{chen2018restricted}
C.~Chen.
\newblock Restricted families of projections and random subspaces.
\newblock {\em Real Anal. Exchange}, 43(2):347--358, 2018.

\bibitem{demeter2020small}
C.~Demeter, L.~Guth, and H.~Wang.
\newblock Small cap decouplings.
\newblock {\em Geom. Funct. Anal.}, 30(4):989--1062, 2020.

\bibitem{duzhang}
X.~Du and R.~Zhang.
\newblock Sharp {$L^2$} estimates of the {S}chr\"{o}dinger maximal function in higher dimensions.
\newblock {\em Ann. of Math. (2)}, 189(3):837--861, 2019.

\bibitem{fassler2014restricted}
K.~F{\"a}ssler and T.~Orponen.
\newblock On restricted families of projections in $\mathbb{R}^3$.
\newblock {\em Proc. Lond. Math. Soc.}, 109(2):353--381, 2014.

\bibitem{fu2021sharp}
Y.~Fu, L.~Guth, and D.~Maldague.
\newblock Sharp superlevel set estimates for small cap decouplings of the parabola.
\newblock {\em arXiv preprint arXiv:2107.13139}, 2021.

\bibitem{gan2022square}
S.~Gan and S.~Wu.
\newblock Square function estimates for conical regions.
\newblock {\em arXiv preprint arXiv:2203.12155}, 2022.

\bibitem{guo2021short}
S.~Guo, Z.~K. Li, P.-L. Yung, and P.~Zorin-Kranich.
\newblock A short proof of $\ell^2$ decoupling for the moment curve.
\newblock {\em Amer. J. Math.}, 143(6):1983--1998, 2021.

\bibitem{GIOW}
L.~Guth, A.~Iosevich, Y.~Ou, and H.~Wang.
\newblock On {F}alconer's distance set problem in the plane.
\newblock {\em Invent. Math.}, 219(3):779--830, 2020.

\bibitem{guth2020improved}
L.~Guth, D.~Maldague, and H.~Wang.
\newblock Improved decoupling for the parabola.
\newblock {\em arXiv preprint arXiv:2009.07953}, 2020.

\bibitem{guth2019incidence}
L.~Guth, N.~Solomon, and H.~Wang.
\newblock Incidence estimates for well spaced tubes.
\newblock {\em Geom. Funct. Anal.}, 29(6):1844--1863, 2019.

\bibitem{guth2020sharp}
L.~Guth, H.~Wang, and R.~Zhang.
\newblock A sharp square function estimate for the cone in $\mathbb{R}^3$.
\newblock {\em Ann. of Math.}, 192(2):551--581, 2020.

\bibitem{harris2019improved}
T.~L.~J. Harris.
\newblock Improved bounds for restricted projection families via weighted fourier restriction.
\newblock {\em arXiv preprint arXiv:1911.00615}, 2019.

\bibitem{harris2021restricted}
T.~L.~J. Harris.
\newblock Restricted families of projections onto planes: the general case of nonvanishing geodesic curvature.
\newblock {\em arXiv preprint arXiv:2107.14701}, 2021.

\bibitem{jarvenpaa2014hausdorff}
E.~J{\"a}rvenp{\"a}{\"a}, M.~J{\"a}rvenp{\"a}{\"a}, and T.~Keleti.
\newblock Hausdorff dimension and non-degenerate families of projections.
\newblock {\em J. Geom. Anal.}, 24(4):2020--2034, 2014.

\bibitem{jarvenpaa2008one}
E.~J{\"a}rvenp{\"a}{\"a}, M.~J{\"a}rvenp{\"a}{\"a}, F.~Ledrappier, and M.~Leikas.
\newblock One-dimensional families of projections.
\newblock {\em Nonlinearity}, 21(3):453, 2008.

\bibitem{kaenmaki2017marstrand}
A.~K{\"a}enm{\"a}ki, T.~Orponen, and L.~Venieri.
\newblock A {M}arstrand-type restricted projection theorem in $\mathbb{R}^3$.
\newblock {\em arXiv preprint arXiv:1708.04859}, 2017.

\bibitem{marstrand1954some}
J.~M. Marstrand.
\newblock Some fundamental geometrical properties of plane sets of fractional dimensions.
\newblock {\em Proc. Lond. Math. Soc.}, 3(1):257--302, 1954.

\bibitem{mattila1975hausdorff}
P.~Mattila.
\newblock Hausdorff dimension, orthogonal projections and intersections with planes.
\newblock {\em Ann. Acad. Sci. Fenn. Ser. AI Math}, 1(2):227--244, 1975.

\bibitem{mattila}
P.~Mattila.
\newblock {\em Geometry of sets and measures in {E}uclidean spaces}, volume~44 of {\em Cambridge Stud. Adv. Math.}
\newblock Cambridge University Press, Cambridge, 1995.
\newblock Fractals and rectifiability.

\bibitem{oberlin2015application}
D.~Oberlin and R.~Oberlin.
\newblock Application of a fourier restriction theorem to certain families of projections in $\mathbb{R}^3$.
\newblock {\em J. Geom. Anal.}, 25(3):1476--1491, 2015.

\bibitem{orponen2015hausdorff}
T.~Orponen.
\newblock Hausdorff dimension estimates for restricted families of projections in $\mathbb{R}^3$.
\newblock {\em Adv. Math.}, 275:147--183, 2015.

\bibitem{orponen2020improved}
T.~Orponen and L.~Venieri.
\newblock Improved bounds for restricted families of projections to planes in $\mathbb{R}^3$.
\newblock {\em Int. Math. Res. Not. IMRN}, 2020(19):5797--5813, 2020.

\bibitem{MR2288738}
M.~Pramanik and A.~Seeger.
\newblock {$L^p$} regularity of averages over curves and bounds for associated maximal operators.
\newblock {\em Amer. J. Math.}, 129(1):61--103, 2005.

\bibitem{pramanik2022furstenberg}
M.~Pramanik, T.~Yang, and J.~Zahl.
\newblock A {F}urstenberg-type problem for circles, and a {K}aufman-type restricted projection theorem in $\mathbb{R}^3$.
\newblock {\em arXiv preprint arXiv:2207.02259}, 2022.

\end{thebibliography}
%***************************************

\end{document}